\documentclass[%
10pt]{amsart}

\usepackage{ifpdf}

\usepackage{amsmath}
\usepackage{amsthm}
\usepackage{amsfonts}
\usepackage{amssymb}

\usepackage[colorlinks,citecolor=blue,linkcolor=blue]{hyperref}

\usepackage{amscd}	
\usepackage{cancel}	
\usepackage{microtype} 

\usepackage{graphicx}
\usepackage{color}

\usepackage[initials]{amsrefs}

\newcommand{\fig}[2]{
\IfFileExists{#1.pdf_tex}{
  \def\svgwidth{#2}\input{#1.pdf_tex}
}{
  \frame{Missing figure ``#1.pdf\_tex''}
  \message{LaTeX Warning: Missing figure ``#1.pdf\_tex'' on input line \the\inputlineno}
}
}

\newcommand{\dimension}{n}

\newcommand{\dist}{\operatorname{dist}}

\newcommand{\argmin}{\operatorname*{arg\,min}}

\newcommand{\divo}{\operatorname{div}}

\newcommand{\dx}{\;dx}

\newcommand{\diff}[1]{\;d{#1}}

\newcommand{\abs}[1]{\left|#1\right|}
\newcommand{\pth}[1]{\left(#1\right)}

\newcommand{\set}[1]{{\left\{#1\right\}}}

\newcommand{\at}[2]{{{\left.{#1}\right|}_{#2}}}
\newcommand{\ang}[1]{{\left\langle#1\right\rangle}}
\newcommand{\norm}[1]{\left\|#1\right\|}

\newcommand{\e}{\ensuremath{\varepsilon}}

\newcommand{\R}{\ensuremath{\mathbb{R}}}

\newcommand{\Rd}{\ensuremath{{\mathbb{R}^{\dimension}}}}
\newcommand{\Rn}{\Rd}
\newcommand{\Z}{\ensuremath{\mathbb{Z}}}
\newcommand{\N}{\ensuremath{\mathbb{N}}}

\definecolor{grey}{rgb}{0.6,0.6,0.6}

\numberwithin{equation}{section}

\swapnumbers	

\theoremstyle{definition}
\newtheorem{example}[subsection]{Example}

\newtheoremstyle{remarkstyle}
  {3pt}
  {3pt}
  {\small}
  {}
  {\bfseries}
  {.}
  { }
  {}

\theoremstyle{remarkstyle}

\usepackage{caption}
\usepackage{subcaption}
\hypersetup{final}
\newcommand{\signdist}{\operatorname{signdist}}
\newcommand{\sign}{\operatorname{sign}}
\newcommand{\shrink}{\operatorname{shrink}}
\newcommand{\Wulff}{\mathcal{W}}

\newcommand{\compl}{{\textsf c}}

\title[Solutions of the crystalline mean curvature flow]{On the self-similar solutions of the crystalline mean curvature flow in three dimensions}

\author[N. Po\v{z}\'{a}r]{Norbert Po\v{z}\'{a}r}
\address[N. Po\v{z}\'{a}r]{Faculty of Mathematics and Physics, Institute of Science and Engineering, Kanazawa University,
Kakuma town, Kanazawa, Ishikawa 920-1192, Japan.}
\email{npozar@se.kanazawa-u.ac.jp}

\date{\today}

\thanks{This work was partially supported by JSPS KAKENHI Grant No. 26800068 (Wakate B) and No. 18K13440 (Wakate).}
\subjclass[2000]{53C44, 35K93, 65M99}
\keywords{crystalline mean curvature flow, level set method, self-similar solutions}

\begin{document}

\maketitle
\begin{abstract}
  We present two types of self-similar shrinking solutions of positive genus for the crystalline mean curvature flow in three dimensions analogous to the solutions known for the standard mean curvature flow.  We use them to test a numerical implementation of a level set algorithm for the crystalline mean curvature flow in three dimensions based on the minimizing movements scheme of A.~Chambolle, \textit{Interfaces Free Bound.~6}~(2004).  We implement a finite element method discretization that seems to improve the handling of edges in three dimensions compared to the standard finite difference method and illustrate its behavior on a few examples.
\end{abstract}

\section{Crystalline mean curvature flow}
The understanding of the evolution of small crystals has been a challenging problem of material science
  and mathematical modeling. In this regime, the evolution seems to be governed by the surface
  energy, whose effects are usually modeled by mean curvature terms. Due to the lattice structure of a
  typical crystal, the surface energy density is
  anisotropic. In fact, it is postulated that the optimal shape (Wulff shape) is a convex polytope and such anisotropies are called crystalline. This causes difficulties for the definition of an
  anisotropic (crystalline) mean curvature and a suitable notion of solutions of the resulting surface evolution
  problem, and makes the development of an efficient numerical method challenging.

  The crystalline mean curvature was introduced independently by S.~B.~Angenent and M.~E.~Gurtin \cite{AG89} and
  J.~E.~Taylor \cite{T91} to model the growth of small crystals, see also \cite{Gurtin,B10}.
The surfaces of solid and liquid bodies have a surface energy, which is usually
expressed as the surface integral of a surface energy density $\sigma: \mathcal{S}^{n-1} \to
(0,\infty)$ over the boundary of a set $E \subset \Rn$, representing the body,
\begin{align*}
  \mathcal{F}(E) := \int_{\partial E} \sigma(\nu) \;dS,
\end{align*}
where $\nu: \partial E \to \mathcal{S}^{n-1}$ is the unit outer normal of $E$. Here $n$ is the
dimension, usually $2$ or $3$ in applications.
For many materials, especially liquids, $\sigma$ is given by the surface tension coefficient and is therefore constant on the unit sphere
$\mathcal{S}^{n-1}$. This surface energy is a manifestation of the fact that the atoms or molecules
forming the body have a smaller interaction energy when surrounded by the particles of the same
kind. Liquids, typically, do not have any preferred direction in the distribution of particles, and
therefore the surface energy density is isotropic.

\begin{figure}
  \centering
  \fig{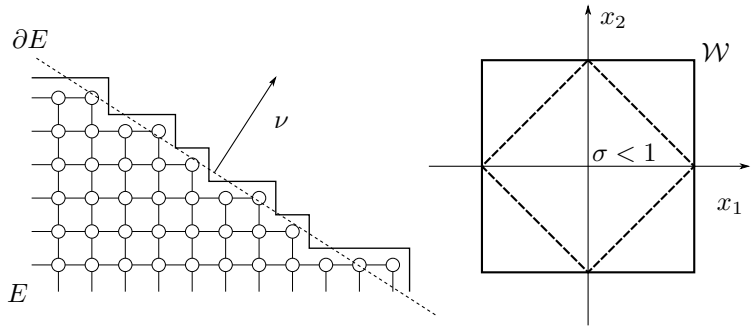}{4in}
  \caption{Left: Crystal of atoms on a regular lattice. Right: The $1$-level set of the associated surface energy density $\sigma$ and its
  Wulff shape $\Wulff$.}
  \label{fig:crystal}
\end{figure}

The situation is quite different for crystals. Let us give a simplified illustration. If we suppose
that the atoms are distributed along the
regular square lattice in two dimensions, every atom inside the body has exactly four neighbors
with which it creates chemical bonds, see Figure~\ref{fig:crystal}. On the surface, however, some of these bonds are broken where
a neighbor is missing and the  surface energy is proportional to the number of the broken bonds.
This number is given in terms of the taxicab ($\ell^1$) length of the surface, not the usual Euclidean
length. In particular, in this case $\sigma(\nu) \sim \norm{\nu}_1 := |\nu_1| + |\nu_2|$, where
$\nu$ is the \emph{macroscopic} unit outer normal. This is
the basic motivation for the introduction of nondifferentiable surface energy densities.

For convenience, we will assume that $\sigma$ is positively one-homogeneously extended from
$\mathcal{S}^{n-1}$ to $\Rn$ as
\begin{align}
  \label{one-homogeneous-extension}
  \sigma(p) = |p| \sigma(\tfrac p{|p|}), \qquad p \in \Rn \setminus \set{0},
\end{align}
$\sigma(0) = 0$,
where $|p| := (\sum_i |p_i|^2)^{1/2}$ is the usual Euclidean norm.

If $\sigma$ is convex on $\Rn$ and $\sigma(p) > 0$ for all $p \in \Rn \setminus \set{0}$, we call
it an \emph{anisotropy}. We are in particular interested in anisotropies $\sigma$ that are
piece-wise linear, for instance the $\ell^1$-norm $\sigma(p) = \norm{p}_1 := \sum_i |p_i|$. Such
anisotropies will be called \emph{crystalline anisotropies}.

The optimal shape of the crystal, that is, the shape with minimal surface energy for a given
volume, is a translation and scaling of the \emph{Wulff shape}
\begin{align*}
  \Wulff := \set{x : x \cdot p \leq \sigma(p), \ p \in \Rn},
\end{align*}
see \cite{Taylor_Symposia}.

The evolution $\{E_t\}_{t\geq 0}$ of a body driven by the dissipation of the surface energy then leads to a formal
gradient flow
\begin{align}
  \label{velocity-law}
  V = \beta(\nu)(\kappa_\sigma + f) \qquad \text{on } \partial E_t,
\end{align}
where $V$ is the normal velocity of the surface $\partial E_t$, $f = f(x,t)$ is a given external
force, and $\beta: \mathcal{S}^{n-1} \to (0, \infty)$ is a \emph{mobility}.
Finally,
$-\kappa_\sigma$ is the first variation of the surface energy $\mathcal{F}$ at $E_t$.
$\kappa_\sigma$ is usually called the anisotropic mean curvature of the surface.

If $\sigma \in C^2(\Rn \setminus \set{0})$ and $\set{p : \sigma(p) < 1}$ is strictly convex, then
it is well-known \cite{G06} that the anisotropic mean curvature can be evaluated as
\begin{align*}
  \kappa_\sigma = -\divo_{\partial E} (\nabla \sigma(\nu)),
\end{align*}
where $\divo_{\partial E}$ is the surface divergence on $\partial E$.

When $\sigma$ is crystalline, the situation is significantly more complicated. In particular, even
if $\partial E$ is smooth,
$\nabla \sigma(\nu)$ might be discontinuous (or not even defined) on parts of the surface.
Therefore $\kappa_\sigma$ might not be defined, or might not be a function.

Instead, following for example \cite{B10}, we define the subdifferential of $\sigma$ as
\begin{align*}
  \partial \sigma(p) := \set{\xi \in \Rn: \sigma(p + h) - \sigma(p) \geq \xi \cdot h,\ h \in \Rn},
\end{align*}
where $\xi \cdot h$ is the usual inner product on $\Rn$. Note that $\partial \sigma(p)$ is a
nonempty compact convex subset of $\Rn$.
We replace $\nabla \sigma(\nu)$ by a vector field $z : \partial E \to \Rn$, usually called a Cahn-Hoffman
vector field, that is a selection of $\partial \sigma(\nu(x))$ on $\partial E$, that is, $z(x) \in
\partial \sigma(\nu(x))$, $x \in \partial E$.
However, now there are multiple choices of $z$ which potentially lead to different values of
$\kappa_\sigma = - \divo_{\partial E} z$. It turns out that a reasonable choice is a vector field $z_{\rm min}$ that
minimizes $\norm{-\divo_{\partial E} z + f}_{L^2(\partial E)}$. The \emph{crystalline (mean) curvature} is then defined as
\begin{align*}
  \kappa_\sigma := - \divo_{\partial E}(z_{\rm min}).
\end{align*}
Such a choice is motivated by the standard theory of monotone operators due to Y.~K{\=o}mura and
H.~Br\'{e}zis \cite{Komura67,Br71}.
Furthermore, since the Euler-Lagrange equation of the minimization problem is $\nabla
(-\divo_{\partial E} z + f) = 0$, this choice yields $\kappa_\sigma + f$ that is constant, if possible, on flat parts, or
facets, of the
crystal parallel to the flat parts of the Wulff shape $\Wulff$. Therefore facets are usually preserved during the evolution, as expected. However,
$\kappa_\sigma$ might be even discontinuous on facets, and then facet breaking or bending occurs,
see
\cite{BNP99} and Figure~\ref{fig:breaking}, Figure~\ref{fig:bending-L}.
This poses a serious difficulty for introducing a suitable notion of solutions for this problem.
Since $\kappa_\sigma$ is itself given as a solution of a minimization problem, it is in general
difficult to evaluate it, except in special circumstances. Moreover, $\kappa_\sigma$ is a nonlocal
quantity on the facets of the crystal as the following example shows.

\begin{example}
  \label{ex:shrinking-cube}
  Consider the cubic anisotropy $\sigma(p) = \norm{p}_1 := \sum_{i=1}^n |p_i|$ and suppose that
  the initial shape is the cube centered at $0$ with side-length $L_0 > 0$, $\beta \equiv 1$. Let us try to find
  $\set{E_t}_{t \geq 0}$. It is not
  difficult to see that for a cube $Q_L = (-\frac L2, \frac L2)^n$, $L > 0$, the vector field $z(x) =
  \frac{x}{\sigma^\circ(x)}$ is a Cahn-Hoffman vector field on $\partial Q_L$, where $\sigma^\circ(x) := \sup
  \set{x \cdot p: \sigma(p) \leq 1}$. Here $\sigma^\circ$ is the convex polar of $\sigma$ and it is the
  dual norm of $\sigma$, $\sigma^\circ(x) = \norm{x}_\infty := \max_{1 \leq i \leq n} |x_i|$, see \cite{Rockafellar}. In
  particular, $\sigma^\circ(x) = \frac L2$ on $\partial Q_L$. Therefore $\divo_{\partial Q_L} z = \frac
  2L (n-1)$. Since it is a constant on the facets, $z$ actually minimizes $\norm{-\divo
  z}_{L^2(\partial Q_L)}$ among all Cahn-Hoffman vector fields, and therefore $\kappa_\sigma = -
  \frac 2L(n-1)$ if $f \equiv const$.  We deduce that, if $\frac 2{L_0}(n-1) \geq f \equiv const$
  so that solution is shrinking, the solution of \eqref{velocity-law} is $\set{E_t}_{t\geq 0}$,
  $E_t = Q_{L(t)}$, where $L(0) = L_0$ and $L' = -\frac 2L(n-1) + f$. When $f = 0$, the unique
  solution is $L(t) = \sqrt{L_0^2 - 4(n-1)t}$. At the \emph{extinction} time $t_1 := \frac{L_0^2}{4(n-1)}$ the cube
  vanishes.

  If the forcing term $f$ is strong enough, the crystal will grow. However, it will only stay a
  cube as long as the velocity of corners is less than $f$ since $\beta \equiv 1$. If the speed of the corners is bigger,
  the corners will round up, as can be easily seen by the comparison principle. See for
  example~\cite{GGR15} and references therein.
\end{example}

\subsubsection*{Solutions of the crystalline mean curvature flow}
Introducing a notion of solutions for \eqref{velocity-law} with the crystalline anisotropy have been
a challenging problem. In two dimensions, if $f$ is constant on facets, the situation is somewhat simpler since $\kappa_\sigma$
is constant on facets of the crystal parallel to the facets of the Wulff shape $\Wulff$. Therefore
if the initial shape is a polygon with edges parallel to edges of  the Wulff shape, the facets will
move without breaking or bending. Their evolution can be tracked by the crystalline algorithm
\cite{T91}, which also
yields efficient numerical methods. However, these methods cannot treat evolutions that are not
strictly faceted.
For fully general situations, the level set method was successfully used to introduce a notion of viscosity solutions to \eqref{velocity-law} in two dimensions by M.-H.~Giga and Y.~Giga \cite{GG98ARMA,GGGakuto,GG01}, and a numerical algorithm was developed by A.~Chambolle~\cite{Chambolle} and further extended to the crystalline case by A.~Oberman, S.~Osher, R.~Takei and R.~Tsai~\cite{OOTT}.

In three dimensions, the situation is significantly more complicated by the possible bending or
breaking of facets.
There is an extensive number of publications that is beyond the scope of this paper, for
instance
\cite{BNP99,BN00,BNP01a,BCCN06}, and \cite{B10} for an introduction to the topic and references.
Recently, A.~Chambolle, M.~Morini and M.~Ponsiglione \cite{CMP}, introduced a well-posed
notion of solutions for the particular velocity law $V = \sigma(\nu) \kappa_\sigma$. In a
subsequent paper with
M.~Novaga \cite{CMNP}, they generalized the theory to $V = \beta(\nu) (\kappa_\sigma + f)$, where
$\beta$
is an anisotropy and $f = f(x,t)$ is a Lipschitz continuous function. They define solutions using the signed distance function to the evolving set, which is required to be a solution of a certain PDE in a sense of distributions, and prove the  existence of such a solution using the minimizing movements algorithm.  Independently,
Y.~Giga and
the author introduced a well-posed notion of viscosity solutions for the level set
formulation of \eqref{velocity-law} in the fully general form $V = F(\nu, \kappa_\sigma + f)$ where
the continuous nonlinearity $F$ is nondecreasing in the second variable, but
with constant $f$ and for bounded crystals \cite{GP16,GP_CPAM}.
See also Section~\ref{sec:level-set}.
It was shown in \cite{CMNP} that both of these notions coincide whenever they both apply.

As for the available numerical results, published results concerning the \emph{purely} crystalline
anisotropy so far seem to only treat the two dimensional
evolution. However, the algorithm proposed in \cite{Chambolle,OOTT} generalizes naturally to three
dimensions, and can easily accommodate a general external force as explained in
Section~\ref{sec:algorithm}. We present some of the results of
this implementation below.
Let us also mention the three dimensional results of J.~W.~Barrett, H.~Garcke and R.~N\"urnberg~\cite{BGN_NM,BGN_IFB,BGN}, who develop a
parametric finite element method for the anisotropic mean curvature flow and
apply it to the Stefan problem with Gibbs-Thomson law that features an
almost-crystalline, but still smooth, anisotropic curvature. This
method does not seem to be able to handle topological changes. For a flow with topological changes their phase field method \cite{BGN_ZAMM,BGN_IMA} is available. See also
\cite{GarckeSurvey} for a survey of numerical approaches.

\subsubsection*{Self-similar solutions}
The study of self-similar solutions of the classical mean curvature flow has been important for the understanding of the singularities of the flow. In two dimensions, it is known that any simple initial closed curve will become a boundary of a convex set in a finite time and therefore the only compact embedded self-similar solution is the circle \cite{GageHamilton,Grayson}.

In three dimension the situation is more interesting. It is known that the sphere is the only convex self-similar solution \cite{Huisken}. The first embedded nonconvex self-similar solution was the ``shrinking doughnut'' solution constructed by S.~B.~Angenent~\cite{Angenent}. This solution can be used to rigorously show the neck-pinching singularity starting with a dumbbell-shaped initial data. Self-similar solutions of higher genus were discovered numerically by D.~L.~Chopp \cite{Chopp}, but to the author's knowledge their existence have not been proven rigorously. The construction of the first higher-genus embedded self-similar solution was done by X.~H.~Nguyen \cite{Nguyen_TransAMS,Nguyen_ADE,Nguyen_Duke}, but this solution is different from the one found by Chopp in \cite{Chopp}.

The behavior is more complex in the crystalline case. There are convex self-similar solutions other than the Wulff shape even in two dimensions. For example, any axes-aligned rectangle will generate a self-similar solution of the crystalline curvature flow with the $\ell^1$ anisotropy. For ``non-rectangular'' even anisotropies, the Wulff shapes in two dimensions are stable \cite{Stancu}. The stability of the Wulff shape solution in three dimensions was considered in \cite{NP_MMMAS}, with more examples of non-Wulff convex solutions. There are also examples of nonconvex self-similar solutions in two dimensions \cite{IUYY}.
For a construction of self-similar solutions in a sector see \cite{GGH}.
The self-similar solutions of positive genus constructed below appear to be new.

\subsection*{Main results and the outline}
We construct the crystalline shrinking doughnut in Section~\ref{sec:shrinking-doughnut} and the crystalline sponge-like solution in Section~\ref{sec:sponge}. In Section~\ref{sec:algorithm} we present an implementation of a numerical algorithm for the crystalline mean curvature flow and test its accuracy using the self-similar solutions.

\section{Self-similar solutions with positive genus}
\label{sec:self-similar}

It is known that the solution of $V = \sigma(\nu) \kappa_\sigma$ with initial data given by the
Wulff shape $\Wulff$ of the anisotropy $\sigma$ is self-similar up to the vanishing time, see
Example~\ref{ex:shrinking-cube}.
In this section we construct self-similar solutions of the anisotropic or crystalline mean curvature
flow with positive genus. One is a torus-like solution, analogous to the solution constructed by S.~B.~Angenent
for the isotropic mean curvature flow \cite{Angenent}, and the other is a sponge-like solution of genus 5, similar to the
solution for the isotropic mean curvature flow discovered numerically by D.~L.~Chopp \cite{Chopp}.
Such solutions are useful for understanding possible singularities of the flow, for instance
showing that a neck-pinching occurs \cite{Angenent}. Furthermore, the solutions can be constructed explicitly for certain
anisotropies and may be therefore useful for testing numerical methods, as we will do in
Section~\ref{sec:numerical-results}.

\subsection{Notion of solutions using the level set method}
\label{sec:level-set}

We need to first clarify what we mean by a solution of the crystalline mean curvature flow \eqref{velocity-law}.
The level set method for the mean curvature flow was introduced and developed in \cite{OS,CGG,ES}.
The basic idea is to introduce an auxiliary function $u: \Rn \times [0, \infty)$, whose
evolution of every level set $\{\set{x \in \Rn: u(x, t) < c}\}_{t \geq 0}$, $c \in \R$, satisfies the velocity law
\eqref{velocity-law}. It is easy to see \cite{G06} that in this case
\begin{align*}
  V = -\frac{u_t}{|\nabla u|}, \qquad \nu = \frac{\nabla u}{|\nabla u|}, \qquad \text{and} \qquad
  \kappa_\sigma = -\divo (\nabla \sigma(\nabla u)).
\end{align*}
Therefore $u$ formally satisfies the equation
\begin{align}
  \label{level-set-eq}
  -\frac{u_t}{|\nabla u|} = \beta\pth{\frac{\nabla u}{|\nabla u|}}(- \divo (\nabla \sigma(\nabla
  u)) + f)
  \qquad \text{in } \Rn \times 0.
\end{align}

If $\sigma$ is a crystalline anisotropy, $\nabla \sigma$ might be discontinuous and therefore the
differential operator on the right-hand side is very singular. In fact, it is a nonlocal operator
on the flat parts of the surface of the crystal parallel to the flat parts of the Wulff shape.
Therefore this equation does not fit within the classical framework of viscosity solutions for
geometric equations \cite{CGG,ES}. The extension of the viscosity theory to \eqref{level-set-eq}
had been a challenging open problem. In one dimension, which also covers two-dimensional crystals,
the theory was developed by M.-H.~Giga, Y.~Giga, P.~Rybka and others \cite{GG98ARMA,GG01,GGR}. Y.~Giga and the
author recently introduced a new notion of viscosity solutions for \eqref{level-set-eq} with $f$
independent of the space variable that
applies to the crystalline anisotropy \cite{GP16,GP_CPAM}. This notion is well-posed for bounded
crystals and stable with
respect to a regularization of the anisotropy, that is, with respect to the approximation of the
crystalline curvature by smooth anisotropic curvatures.
The main idea is to restrict the space of test functions to \emph{faceted} test functions for which it is possible define the operator $\divo [\nabla \sigma(\nabla u)]$ as the
divergence of a minimizing Cahn-Hoffman vector field $z \in \partial \sigma(\nabla u)$ as explained above. The solutions constructed below are viscosity solutions of the level set formulations.

\subsection{The shrinking doughnut}
\label{sec:shrinking-doughnut}

In \cite{Angenent}, S.~B.~Angenent showed the existence of a doughnut-like self-similar shrinking
solution for the isotropic mean curvature flow in dimensions $n \geq 3$.

An analogue of this solution can be also constructed in the anisotropic case.
Let us fix a dimension $n \geq 3$.
We consider a ``cylindrical'' anisotropy
\begin{align*}
  \sigma(p) := \tilde \sigma(p') + |p_n|, \qquad p = (p', p_n) \in \Rn,
\end{align*}
where $\tilde\sigma$ is an even anisotropy on $\R^{n-1}$ (not necessarily crystalline), i.e.,
we require that $\tilde\sigma(-p') = \tilde\sigma(p')$.
We define the mobility
\begin{align}
  \label{beta-torus}
  \beta(p) := \tilde\sigma(p') + \mu |p_n|, \qquad p = (p', p_n) \in \Rn,
\end{align}
where
$\mu > 0$ is a given constant. We will see later that we must take $\mu = \frac 12$ for $n = 3$ to get a self-similar solution.
Note that
\begin{align*}
  \sigma^\circ(x) = \max(\tilde\sigma^\circ(x'), |x_n|), \qquad x = (x', x_n) \in \Rn,
\end{align*}
where $\tilde\sigma^\circ(x') := \sup \set{x'\cdot p' : p' \in \R^{n-1},\ \tilde\sigma(p') \leq 1}$ is
the convex polar of $\tilde\sigma$,
see for instance \cite{Rockafellar} for this and other convex analysis results.

\begin{figure}
  \centering
\begin{subfigure}{0.48\textwidth}
  \fig{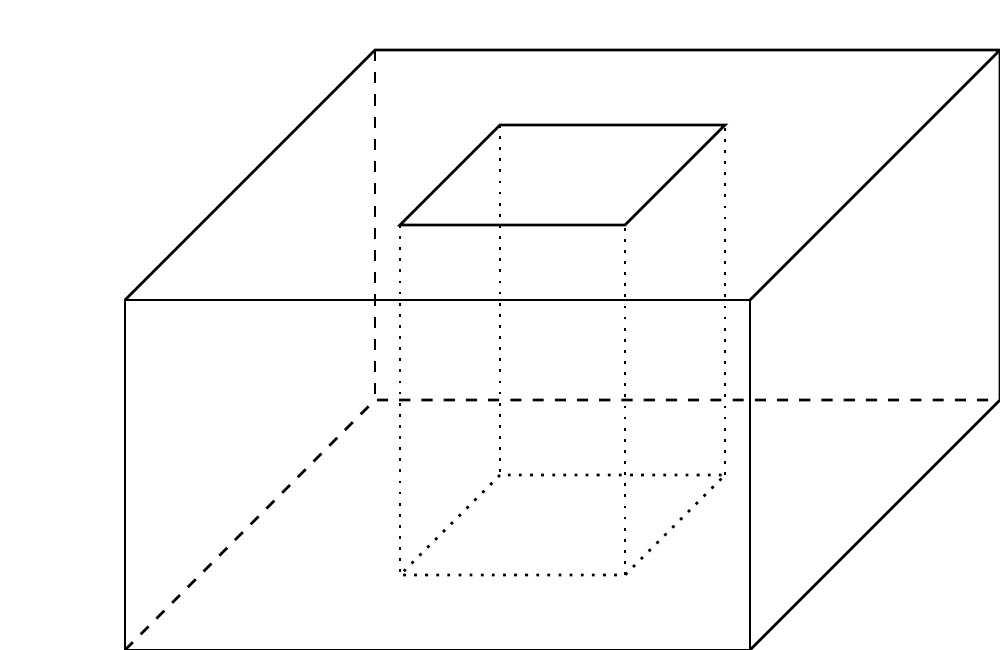}{\textwidth}
\end{subfigure}
\begin{subfigure}{0.48\textwidth}
  \includegraphics[width=\textwidth]{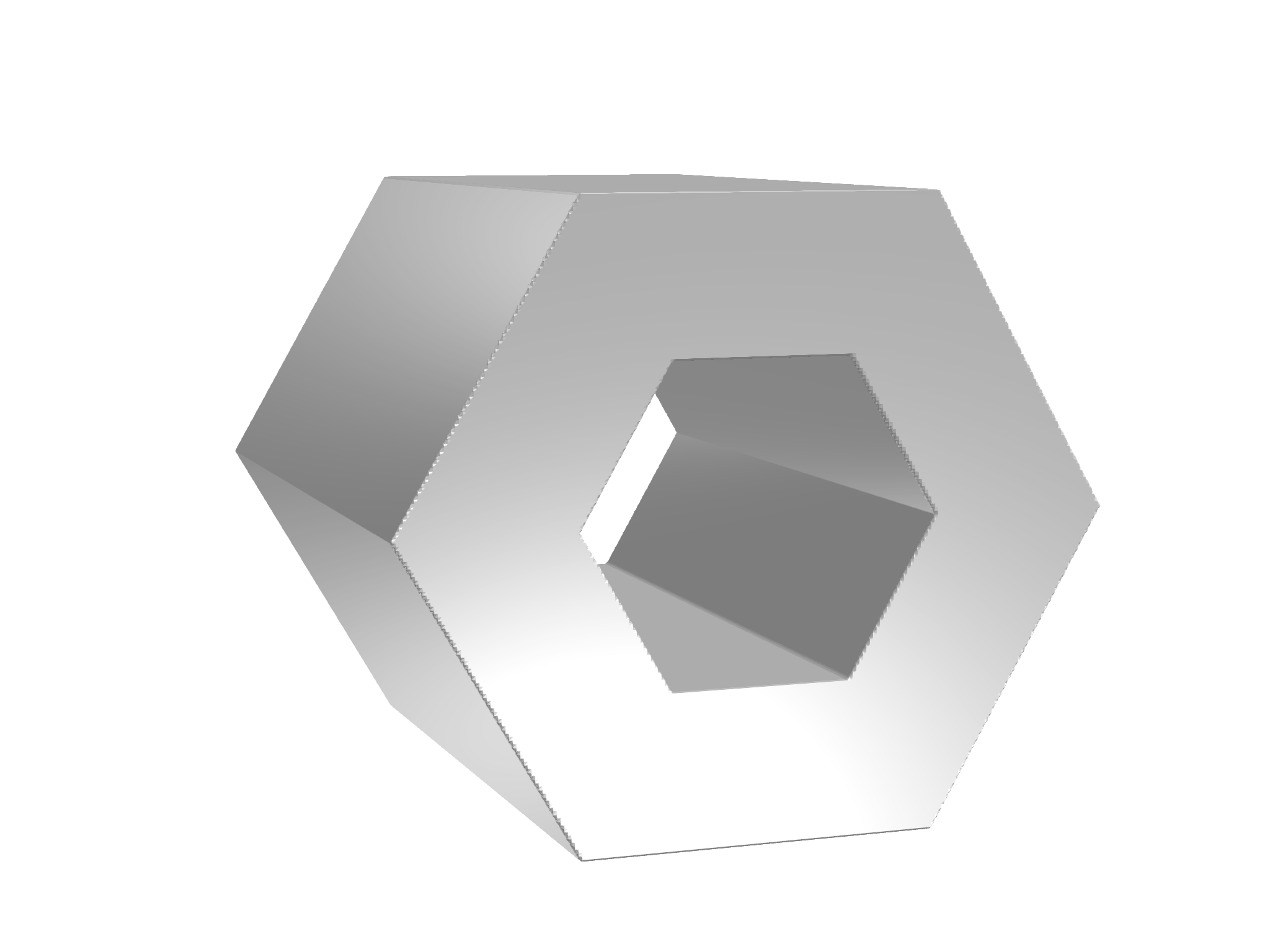}
\end{subfigure}
  \caption{Torus $T_{r, R, h}$ for $n = 3$ and $\tilde\sigma(p') = |p'|_1$ (left), and with hexagonal anisotropy (right).}
  \label{fig:torus}
\end{figure}

We consider the anisotropic mean curvature flow
\begin{align}
  \label{torus-vel}
  V = \beta(\nu) \kappa_\sigma.
\end{align}
Let $R > r > 0$, $h > 0$ be real parameters. Define the set
\begin{align*}
  T_{r,R,h} := \set{x = (x', x_n) \in \R^n: r < \tilde\sigma^\circ(x') < R,\ |x_n| < h}.
\end{align*}
This is a ``torus'' with the hole aligned with the $x_n$-direction, see Figure~\ref{fig:torus}.
We will use this set to construct a self-similar evolution
\begin{align*}
  E_t :=
  \begin{cases}
    T_{r(t), R(t), h(t)}, & 0\leq t < t^*,\\
    \emptyset, & t \geq t^*,
  \end{cases}
\end{align*}
given appropriate functions $r(t)$, $R(t)$ and $h(t)$, where $t^*$ is the extinction time.

Let us first calculate the anisotropic curvature on the surface of $T_{r,R,h}$.
We will split the surface $\Gamma = \partial T_{r, R, h}$ into (o) the outer surface $\Gamma_{\rm
o}$ where
$\tilde\sigma^\circ(x') = R$, (i) the inner surface $\Gamma_{\rm i}$
where $\tilde\sigma^\circ(x') = r$, and (s)
the side facets $\Gamma_{\rm s}$ with $|x_n| = h$.

To compute the anisotropic curvature, we need to find a Cahn-Hoffman vector field $z: \Gamma \to
\Rn$ on the surface
such that
$z(x) \in \partial\sigma(\nu(x))$ for $\mathcal H^{n-1}$-a.e. $x$ with minimal
$\norm{\divo_{\Gamma} z}_{L^2(\Gamma)}$.
It can be shown that $\divo_{\Gamma} z$ minimizes the $L^2$ norm exactly when $\divo_{\Gamma} z$ is
constant on the flat parts of $\Gamma$, see \cite{B10}.

One such $z$ can be given explicitly as
\begin{align}
  \label{torus-z}
  z(x) = \pth{g(\tilde\sigma^\circ(x'))x', \frac{x_n}{h}}, \qquad x = (x', x_n) \in \Gamma,
\end{align}
where
\begin{align*}
  g(s) := a s^{-n + 1} +b,
\end{align*}
with
\begin{align*}
  a = -\frac{R^{n-2} r^{n-2}(R + r)}{R^{n-1} - r^{n-1}},\qquad
  b = \frac{R^{n-2} + r^{n-2}}{R^{n-1} - r^{n-1}}.
\end{align*}
The constants $a$ and $b$ are chosen so that $g(r)r = -1$ and $g(R)R = 1$. Since $r < R$ and
$g(s)s$ is increasing in $s$ as $a < 0$, we conclude that $-1 < g(s) s < 1$ for $r < s < R$.
In particular,
\begin{align}
  \label{g-choice}
  \tilde\sigma^\circ(g(\tilde\sigma^\circ(x'))x') \leq 1 \Leftrightarrow g(\tilde\sigma^\circ(x'))x' \in \partial
  \tilde\sigma(0) \qquad \text{for all } x = (x', x_n) \in \Gamma.
\end{align}
Furthermore,
\begin{align*}
  g'(s) s +(n-1) g(s) = (n-1)b
\end{align*}
and thus
\begin{align}
  \label{gpart-div}
  \begin{aligned}
  \divo_{x'} (g(\tilde\sigma^\circ(x'))x') &= g'(\tilde\sigma^\circ(x')) x' \cdot \nabla
  \tilde\sigma^\circ(x') + (n-1)g(\tilde\sigma^\circ(x'))\\
  &=g'(\tilde\sigma^\circ(x'))\tilde\sigma^\circ(x')  + (n-1)g(\tilde\sigma^\circ(x'))
  =(n-1)b
  \end{aligned}
\end{align}
whenever $\nabla\tilde\sigma^\circ(x')$ exists, where we used that $x' \cdot \nabla \tilde\sigma^\circ(x')
= \tilde\sigma^\circ(x')$, see \cite{Rockafellar}.

Let us check that $z$ is indeed a Cahn-Hoffman vector field on $\Gamma$.
A convenient approach is to write $\Gamma$ as the level set of a Lipschitz continuous function.
Consider
\begin{align}
  \label{psi-level-set}
  \psi(x) = \max(r - \tilde\sigma^\circ(x'), \tilde\sigma^\circ(x') - R, |x_n| - h),
  \qquad x = (x', x_n) \in \Rn.
\end{align}
Clearly $\psi$ is Lipschitz, $\Gamma = \set{x: \psi(x) = 0}$ and $T_{r,R,h} = \set{x: \psi(x) < 0}$. $\nabla \psi$ is an outer normal
vector of $\Gamma$ (with respect to $T_{r, R, h}$) and exists $\mathcal H^{n-1}$-a.e. on $\Gamma$.
Since $\partial \sigma$ is positively zero-homogeneous, it is enough to show that $z(x) \in
\partial \sigma(\nabla \psi(x))$ $\mathcal H^{n-1}$-a.e. $x \in \Gamma$.

Let $\operatorname{ri} \Gamma_j$ denote the relative interior of the surface $\Gamma_j$.
Thus suppose that $x \in \Gamma$ such that $\nabla \psi(x)$ exists.
We have three cases:
\newcommand{\ri}{\operatorname{ri}}
\begin{description}
  \item[$x \in \ri\Gamma_{\rm o}$] In a neighborhood of such a point, we see that $\psi(y) =
    \tilde\sigma(y') - R$, and so in this case we have $\nabla \psi(x) = (\nabla
    \tilde\sigma^\circ(x'), 0)$, from which we deduce
    \begin{align*}
      z(x) = \pth{\frac{x'}{\tilde\sigma^\circ(x')}, \frac{x_n}{h}} \in \partial \tilde\sigma(\nabla
      \tilde\sigma^\circ(x')) \times [-1, 1] = \partial\sigma(\nabla\psi(x)).
    \end{align*}
  \item[$x \in \ri\Gamma_{\rm i}$] This can be handled as the previous case, recalling that
    $\tilde\sigma$ is even and thus $-\frac{x'}{\tilde\sigma^\circ(x')} \in \partial
    \tilde\sigma(-\nabla \tilde\sigma^\circ(x'))$.
  \item[$x \in \ri\Gamma_{\rm s}$] Now we are on the top or the bottom flat facet, $\psi(y) = |y_n| -
    h$ in the neighborhood of this point, and so $\nabla \psi(x) = \pth{0, \frac{x_n}{|x_n|}}$.
    Thus, recalling \eqref{g-choice}, we deduce
    \begin{align*}
      z(x) = \pth{g(\tilde\sigma^\circ(x'))x', \frac{x_n}{|x_n|}} \in \partial \tilde\sigma(0) \times
      \set{\frac{x_n}{|x_n|}} = \partial\sigma(\nabla\psi(x)).
    \end{align*}
\end{description}
We have proved that $z$ defined in \eqref{torus-z} is a Cahn-Hoffman field on $\Gamma$.

Now we show that $\divo_\Gamma z$ is constant on the flat parts of $\Gamma$, and hence it minimizes
$\norm{\divo_\Gamma z}_{L^2(\Gamma)}$.
The surface (tangential) divergence $\divo_\Gamma z(x)$ for $x \in \Gamma$ can be computed easily using the level set
method by constructing an
extension $\bar z$ of $z$ away
from the facet such that $\bar z \in \partial \sigma(\nabla \psi)$ a.e. in a neighborhood of
$x \in \Gamma$ so that $\divo \bar z(x)$ exists. Then we have $\divo_\Gamma z(x) =\divo \bar z(x)$,
see \cite{G06}.

Following
the above verification that $z$ is a Cahn-Hoffman vector field, we can construct the extension
$\bar z$ in a neighborhood of $x \in \Gamma$ in
the following three cases in the following way:
\begin{description}
  \item[$x\in \operatorname{ri} \Gamma_{\rm o}$] Here we can simply take $\bar z(y) =
    \pth{\frac{y'}{\tilde\sigma^\circ(y')}, \frac{y_n}{h}}$ which yields
    \begin{align*}
      \divo_\Gamma z(x) = \divo \bar z(x) = \frac{n-2}R + \frac1h,
    \end{align*}
    if the divergence exists,
    where we followed the computation in \eqref{gpart-div} with $g(s) = \frac 1s$.
  \item[$x\in \operatorname{ri} \Gamma_{\rm i}$] Similarly to the previous situation, we may take
  \begin{align*}\bar z(y) = \pth{-\frac{y'}{\tilde\sigma^\circ(y')}, \frac{y_n}{h}},\end{align*}
    which yields
    \begin{align*}
      \divo_\Gamma z(x) = \divo \bar z(x) = -\frac{n-2}r + \frac1h
    \end{align*}
    if the divergence exists.
  \item[$x\in \operatorname{ri} \Gamma_{\rm s}$] Here we may take $\bar z(y) =
    \pth{g(\tilde\sigma^\circ(y'))y', \frac{y_n}{|y_n|}}$, yielding
    \begin{align*}
      \divo_\Gamma z(x) = \divo \bar z(x) = (n-1)b
    \end{align*}
    using \eqref{gpart-div}.
\end{description}
Since $\divo_\Gamma z$ is constant on the flat parts of the surface, we deduce that it minimizes the
$L^2$ norm among all Cahn-Hoffman vector fields and therefore we can define $\kappa_\sigma = -
\divo_\Gamma z$.

Let us now find the normal velocity $V$ of the surface $\Gamma_t = \partial T_{r(t), R(t), h(t)}$
where $r(t)$, $R(t)$ and $h(t)$ are some $C^1$ functions satisfying $0 < r(t) < R(t)$, $0 < h(t)$.
It is straightforward to use the level set method. Let $\psi(x,t)$ be the function defined in
\eqref{psi-level-set} with the time dependent parameters $r(t)$,  $R(t)$, $h(t)$ above. Let us fix $t$ and $x \in \Gamma_t$
so that $\psi_t(x,t)$ and $\nabla \psi(x, t)$ are defined. We can again consider three cases
and use the formula $V(x,t) = -\frac{\psi_t}{|\nabla \psi|}(x,t)$ \cite{G06}:
\begin{align}
  \label{torus-V}
  V(x,t) =
  \begin{cases}
    \frac{R'(t)}{|\nabla \tilde\sigma^\circ(x')|}, & x \in \ri \Gamma_{\rm o},\\
    -\frac{r'(t)}{|\nabla \tilde\sigma^\circ(x')|}, & x \in \ri \Gamma_{\rm i},\\
    h'(t), & x \in \ri \Gamma_{\rm s}.
  \end{cases}
\end{align}

To relate $V$ to $\kappa_\sigma$ using the law \eqref{torus-vel}, it is left to compute
$\beta(\nu)$ on $\Gamma$. Using the level set method again, $\nu = \frac{\nabla
\psi}{|\nabla\psi|}$ whenever $\nabla \psi$ exists and is nonzero \cite{G06}, and therefore
\begin{align}
  \label{torus-nu}
  \nu(x) =
  \begin{cases}
    \frac{(\nabla \tilde\sigma^\circ(x'),0)}{|\nabla \tilde\sigma^\circ(x')|}, & x \in \ri \Gamma_{\rm o},\\
    -\frac{(\nabla \tilde\sigma^\circ(x'),0)}{|\nabla \tilde\sigma^\circ(x')|}, & x \in \ri \Gamma_{\rm i},\\
    (0, \frac{x_n}{|x_n|}), & x \in \ri \Gamma_{\rm s}.
  \end{cases}
\end{align}
Recalling $\beta$ given in \eqref{beta-torus} and that $\tilde\sigma(\pm\nabla \tilde\sigma^\circ(x')) = 1$
whenever $\nabla\tilde\sigma^\circ(x')$ exists \cite{Rockafellar}, we deduce
\begin{align}
  \label{torus-beta}
  \beta(\nu(x)) =
  \begin{cases}
    \frac{1}{|\nabla \tilde\sigma^\circ(x')|}, & x \in \ri \Gamma_{\rm o} \cup \ri \Gamma_{\rm i},\\
    \mu, & x \in \ri \Gamma_{\rm s}.
  \end{cases}
\end{align}

Relating the anisotropic curvature $\kappa_\sigma = - \divo_\Gamma z$ with \eqref{torus-V}, \eqref{torus-nu} and \eqref{torus-beta} via \eqref{torus-vel}, we must
have
\begin{align}
  \label{torus-odes}
  \begin{aligned}
  R' &= -\frac{n-2}R - \frac 1h,\\
  r' &= -\frac{n-2}r + \frac 1h,\\
h' &= -(n-1)\mu \frac{R^{n-2} + r^{n-2}}{R^{n-1} - r^{n-1}}.
  \end{aligned}
\end{align}

For the evolution to be self-similar, we therefore need
\begin{align*}
  \pth{\frac rR}' = 0 = \pth{\frac hR}'.
\end{align*}
Thus there are constants $\gamma \in (0,1)$ and $\lambda > 0$ such that
\begin{align*}
  r = \gamma R, \qquad h= \lambda R.
\end{align*}
Plugging this into \eqref{torus-odes} and multiplying each equation by $R$, we get
\begin{subequations}
\begin{align}
  \label{alg-a}
  RR' &= -(n-2)-\frac1\lambda\\
  \label{alg-b}
  &=-\frac{n-2}{\gamma^2} + \frac 1{\lambda\gamma}\\
  \label{alg-c}
  &=-\frac{(n-1)\mu}\lambda\cdot \frac{1+\gamma^{n-2}}{1 - \gamma^{n-1}}.
\end{align}
\end{subequations}

From \eqref{alg-a} and \eqref{alg-b}, we obtain
\begin{align}
  \label{lambda-1}
  \lambda = \frac{\gamma}{(n-2)(1-\gamma)}.
\end{align}
On the other hand, \eqref{alg-a} and \eqref{alg-c} yield
\begin{align}
  \label{lambda-2}
  \lambda = \frac 1{n-2}\pth{(n-1)\mu\frac{1+\gamma^{n-2}}{1 - \gamma^{n-1}}-1}.
\end{align}
Hence combining \eqref{lambda-1} and \eqref{lambda-2} and multiplying both sides by $(n-2)(1 -
\gamma^{n-1})$, we have
\begin{align*}
  1 + \gamma +\cdots + \gamma^{n-2}  = (n-1)\mu (1 + \gamma^{n-2}).
\end{align*}

This equation has a solution $\gamma \in (0,1)$ only for certain $\mu$.
Let us only consider $n = 3$. Then we have
\begin{align*}
  0 = (2\mu - 1) (1 + \gamma),
\end{align*}
which has a solution $\gamma \in (0,1)$ only if $\mu = \frac 12$.
In this case there are infinitely many solutions,
for instance,
\begin{align*}
  \gamma = \tfrac 12,\quad \lambda = 1.
\end{align*}
Then the solution of \eqref{torus-vel} with initial data $E_0 = T_{\frac12, 1,1}$ is
the self-similar evolution
\begin{align*}
  E_t =
  \begin{cases}
    \sqrt{1 - 4t} E_0, & 0 \leq t < t^* = \frac14,\\
    \emptyset, & t^* \leq t.
  \end{cases}
\end{align*}

It is possible to show that this is the unique (open) evolution given by the viscosity solution of
\eqref{torus-vel} in the sense of \cite{GP_CPAM}, but this is beyond the scope of this note. For
more details, see \cite{GGGakuto}.

\subsection{The sponge}
\label{sec:sponge}
\begin{figure}
  \begin{subfigure}{0.24\textwidth}
  \includegraphics[width=\textwidth]{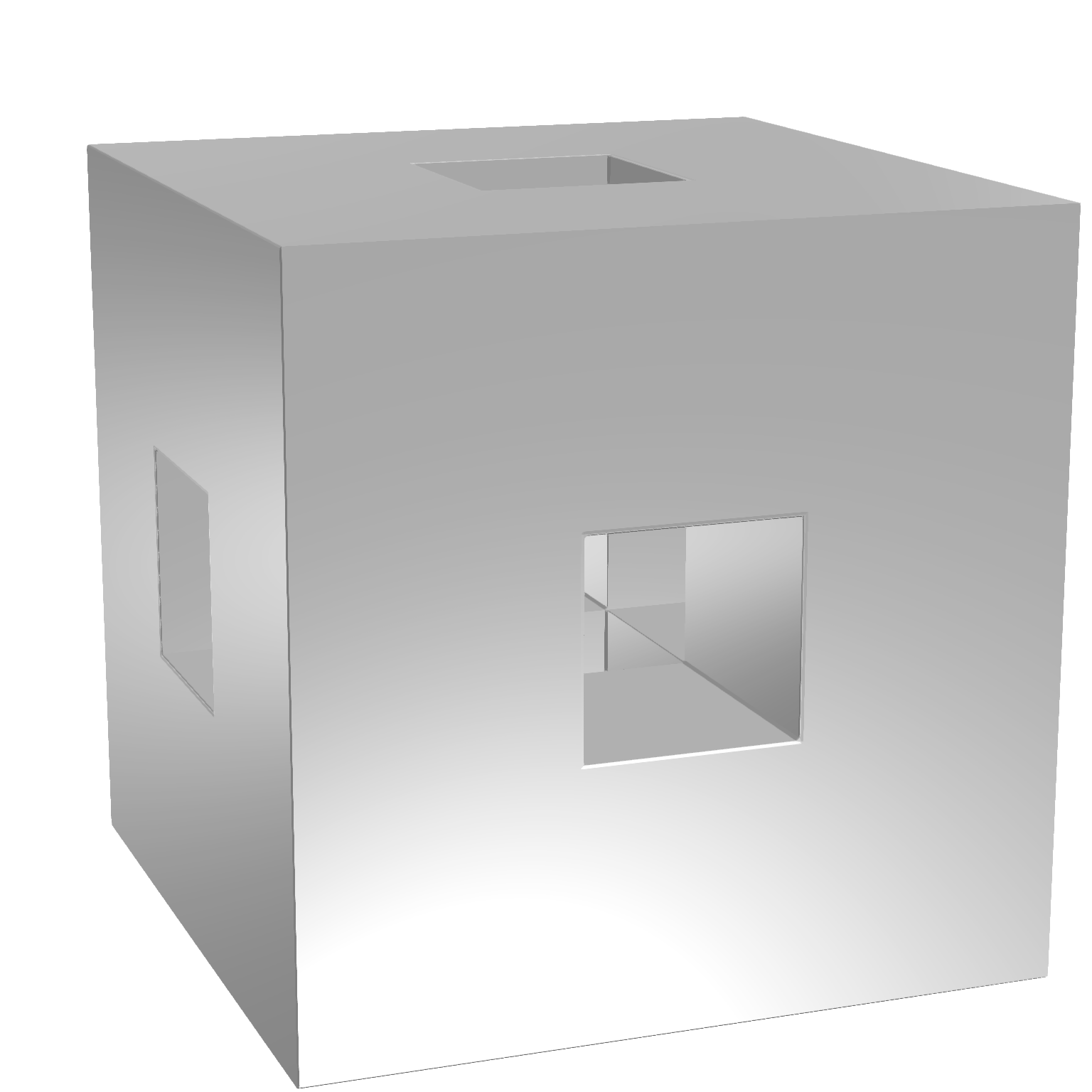}
  \caption{}
  \end{subfigure}
  \begin{subfigure}{0.24\textwidth}
    \includegraphics[width=\textwidth]{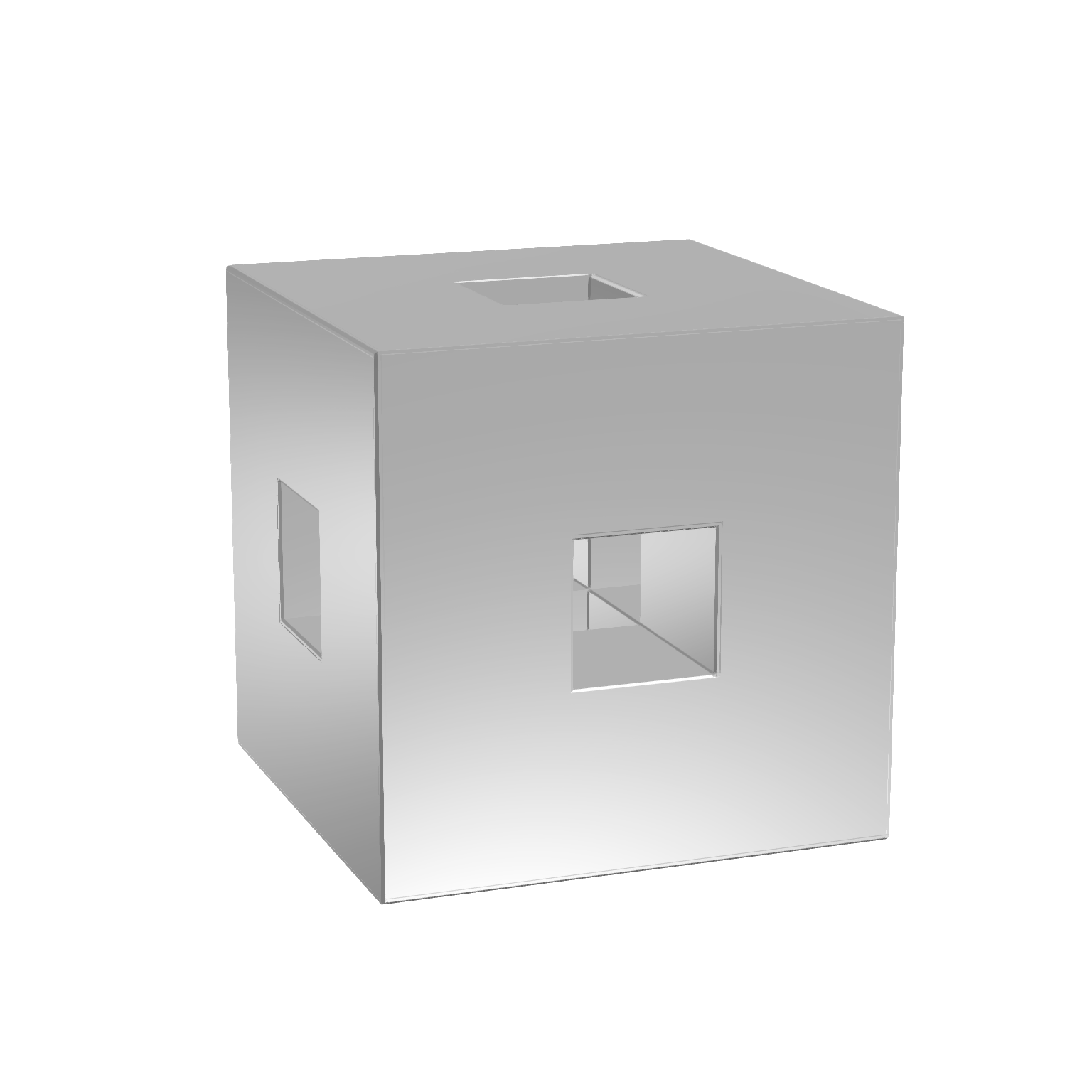}
  \caption{}
  \end{subfigure}
  \begin{subfigure}{0.24\textwidth}
  \includegraphics[width=\textwidth]{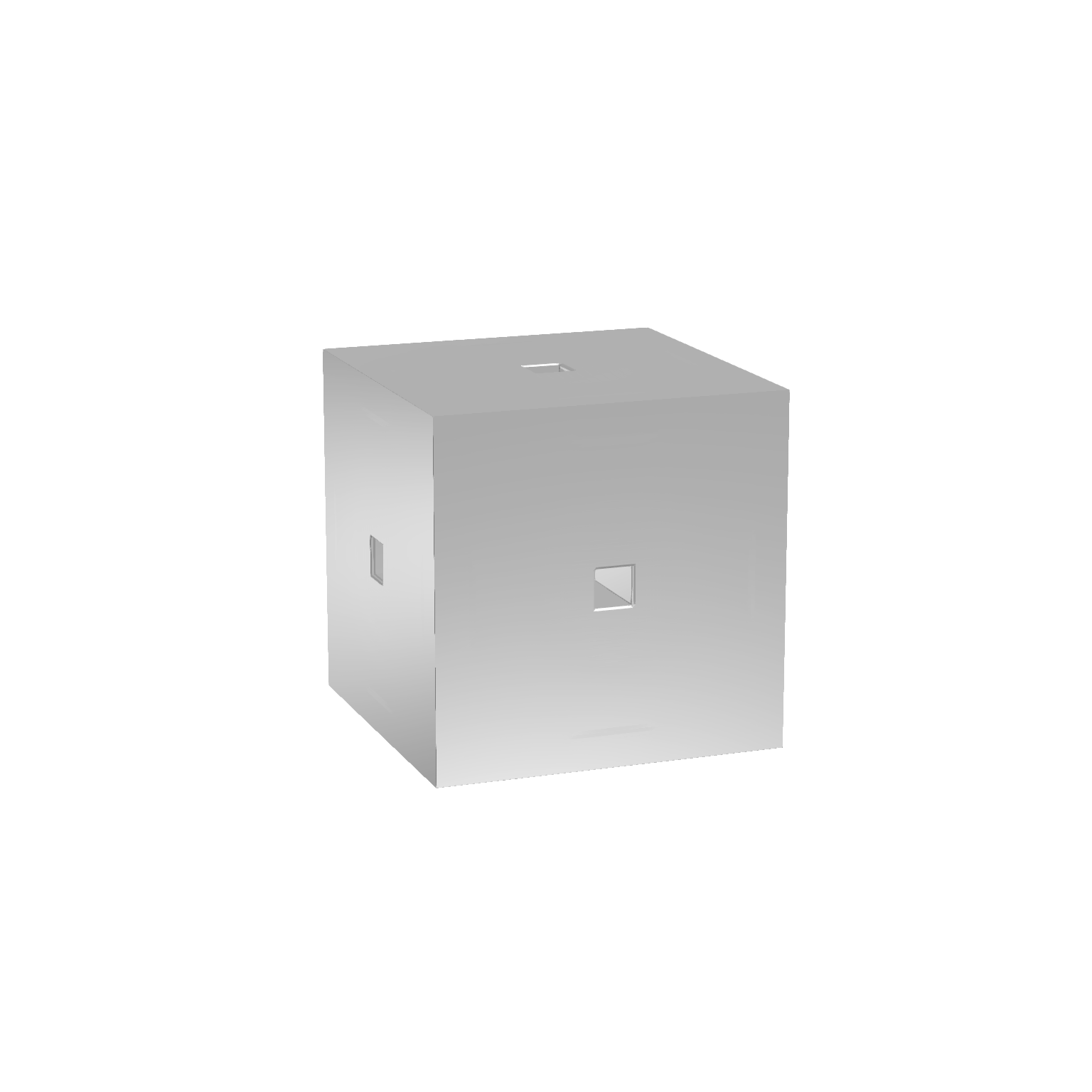}
  \caption{}
  \end{subfigure}
  \begin{subfigure}{0.24\textwidth}
    \includegraphics[width=\textwidth]{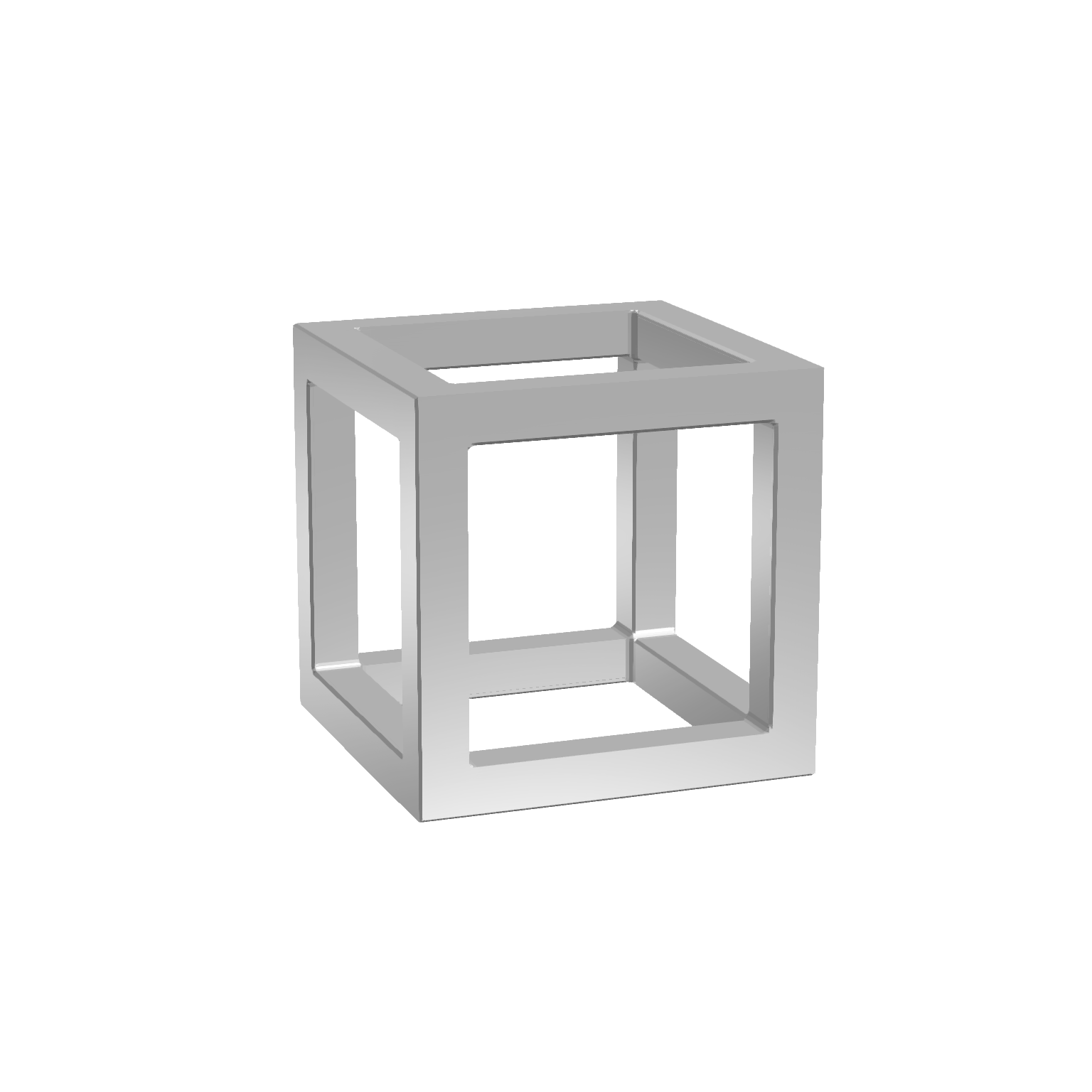}
  \caption{}
  \end{subfigure}
  \caption{
    The shrinking sponge.
    (A) The initial set $E_0$. The outer edge has length $2R_0$ while the inner edge has length $2r_0$.
    (B) Self-similar sponge solution with initial size $0.8$ at time $t = 0.015$
    computed with domain resolution $M = 256$.
    A numerical solution uses Chambolle's algorithm; since the solution
      is unstable, the accumulated numerical errors will cause it
    to eventually diverge from the self-similar solution. (C) Evolution for $\frac{R_0}{r_0} >
  \xi^*$. (D) Evolution for $\frac{R_0}{r_0} < \xi^*$.}
  \label{fig:sponge}
\end{figure}

In this section we construct
a shrinking ``sponge'' in dimension $n = 3$. It is a self similar
solution of the crystalline mean curvature flow with cubic anisotropy $\sigma(p) = |p|_1$ with genus $5$ and resembles the level $1$ Menger sponge. Its shape
is given by a cube centered at the origin from which the neighborhood (in the maximum
($\ell^\infty$) norm) of all
three axes has been removed, see Figure~\ref{fig:sponge}.
It is an analogue of the solution that was discovered numerically by D.~L.~Chopp
\cite{Chopp} for the standard isotropic mean curvature flow.

Let us thus consider the crystalline mean curvature flow
\begin{align}
  \label{sponge-vel}
  V = \sigma(\nu)\kappa_\sigma
\end{align}
where $\sigma(p) = |p|_1$ is the cubic ($\ell^1$) anisotropy, with initial data
\begin{align}
  \label{initial}
  E_0 := S_{r_0, R_0}
\end{align}
for some constants $0 < r_0 < R_0$,
where
\begin{align*}
  S_{r,R} &:=
  \set{x \in \R^3: |x|_\infty < R,\ |x_i| > r \text{ for at least two of $i =
  1,2,3$}}
\end{align*}
for $0 < r < R$,
see Figure~\ref{fig:sponge}.

We claim that the open evolution $\set{E_t}_{t \geq 0}$ solving \eqref{sponge-vel} is of the form
\begin{align*}
  E_t := S_{r(t), R(t)},
\end{align*}
where $r(t)$, $R(t)$ are solutions of a certain ODE system.
Our goal is to find $0 < r_0 < R_0$ so that the evolution $\set{E_t}_{t \geq 0}$ is self-similar.

For $0 < r <R$, due to the symmetries, $S_{r, R}$ has only exactly two types of facets: (o) the outer facet given by a square of side
$2R$ with square of side $2r$ with the same center removed, and (i) the inner facet given by a rectangle
with sides of lengths $2r$ and $R - r$.
The crystalline curvature on each of these facets can be expressed as the ratio of their \emph{signed}
perimeter and their area, that is,
\begin{align*}
  \kappa_{\textrm o} &=  -\frac{8(R + r)}{4R^2 - 4r^2} = - \frac 2{R - r},&
  \kappa_{\textrm i} &= \frac 1r - \frac 2{R-r} = \frac{R - 3r}{r(R-r)}.
\end{align*}
For a detailed explanation of why these are indeed the curvatures, see
Section~\ref{sec:shrinking-doughnut} for an explicit Cahn-Hoffman field construction, or \cite{LMM}.

Following the reasoning in Section~\ref{sec:shrinking-doughnut}, the normal velocity of the outer facet of $\set{E_t}$
is $V = R'(t)$ and the inner facet is $V = -r'(t)$.
Note also that since the facets of $S_{r,R}$ are axes-aligned, we have $\sigma(\nu) = |\nu|_1 = 1$
a.e. on $\partial S_{r,R}$.
Therefore \eqref{sponge-vel} simplifies for $\set{E_t}$ into the system of ODEs
\begin{align}
  \label{sponge-crystalline-algorithm}
  \left\{\begin{aligned}
  R' &= \kappa_{\textrm o} = -\frac 2{R-r},\\
  r' &= -\kappa_{\textrm i} = -\frac {R-3r}{r(R-r)},
  \end{aligned}\right.
  \qquad t > 0,
\end{align}
with initial data $r(0) = r_0$, $R(0) = R_0$.

Let us perform a more detailed analysis of the behavior of the system
\eqref{sponge-crystalline-algorithm}.
The evolution will be self-similar if
\begin{align*}
  \xi(t) := \frac{R(t)}{r(t)} \equiv const.
\end{align*}
Rewriting the system \eqref{sponge-crystalline-algorithm} for $\xi$ and $r$, we have
\begin{align}
  \label{sponge-crystalline-algorithm-xi}
  \left\{\begin{aligned}
      \xi' &= \frac{\xi^2 - 3 \xi - 2}{r^2 (\xi - 1)},\\
      r' &= - \frac{\xi - 3}{r ( \xi -1)},
  \end{aligned}\right.
  \qquad t > 0,
\end{align}
with initial data $\xi_0 = \frac{R_0}{r_0} > 1$, $r(0) = r_0 > 0$.

We set
\begin{align*}
  g(\xi) := \xi^2 - 3 \xi - 2.
\end{align*}
The constant $\xi^* := \frac{3 + \sqrt{17}}2$ is the unique value $> 1$ such that $g(\xi^*) = 0$.
Let us also introduce $Q_R := \set{x \in \R^3: |x|_\infty < R}$

We can classify the behavior of the solution of \eqref{sponge-crystalline-algorithm-xi} according
to the initial ratio $\xi_0 = \frac{R_0}{r_0}$ as follows:
\begin{enumerate}
  \item $1 < \xi_0 <  \xi^*$: Since $\xi' < 0$ and $r(t) \leq R_0$, there exists a time $t^* >
    0$ such that $\lim_{t \to t^*-} \xi(t)  = 1$, that is, $\lim_{t \to t^*-} R(t) - r(t)  = 0$.
    The sponge gets thinner with time and converges to the edges of the cube
    $Q_{R(t^*)}$, and vanishes at $t^*$; Figure~\ref{fig:sponge}(D).
  \item $\xi_0 = \xi^*$: The solution is $\xi(t) = \xi_0 = \xi^*$, $r(t) = \sqrt{r_0^2 - 2ct}$ for $0 < t <
    t^*$, where $c = \frac{\xi^* - 3}{\xi^* - 1} > 0$ and $t^* = \frac{r_0^2}{2c} > 0$. At $t =
    t^*$ the sponge vanishes. The evolution is \emph{self-similar}.
  \item $\xi_0 > \xi^*$: There is time $t^\Box > 0$ such that
    $\lim_{t \to t^\Box-} (\xi(t), r(t), R(t)) = (+\infty, 0, R(t^\Box))$. In this case, the holes
    close up and the sponge becomes a cube at time $t^\Box$; Figure~\ref{fig:sponge}(D). After $t = t^\Box$, $R$ evolves with
    $R' = - \frac{2}{R}$. The cube vanishes at a later time $t^* > t^\Box$.
\end{enumerate}

See Figure~\ref{fig:sponge} for a numerical solution using Chambolle's algorithm discussed in
Section~\ref{sec:algorithm}. Note that (b) is
unstable since $g'(\xi^*) > 0$, and therefore numerically we will observe either (a) or (c).

Let us summarize the solution of \eqref{sponge-vel} in terms of the values $0 < r_0 < R_0$:
\begin{description}
  \item[$1 < \frac{R_0}{r_0} \leq \xi^*$]
    \begin{align*}
       E_t :=
       \begin{cases}
          S_{r(t), R(t)}, & 0 \leq t < t^*,\\
          \emptyset, & t \geq t^*.
       \end{cases}
    \end{align*}
  \item[$\frac{R_0}{r_0} > \xi^*$]
    \begin{align*}
       E_t :=
       \begin{cases}
          S_{r(t), R(t)}, & 0 \leq t < t^\Box,\\
          Q_{R(t)}, & t^\Box \leq t < t^*,\\
          \emptyset, & t \geq t^*.
       \end{cases}
    \end{align*}
\end{description}

As in the shrinking doughnut case by following \cite{GGGakuto}, one can show that this is the unique (open) evolution given by
the viscosity solution of \eqref{sponge-vel} in the sense of \cite{GP_CPAM}.

\section{The numerical algorithm}
\label{sec:algorithm}

An efficient method for the mean curvature flow \eqref{velocity-law}
is based on a minimizing movement formulation due to A.~Chambolle \cite{Chambolle}, which can be efficiently
solved by a split Bregman iteration proposed by \cite{GO,OOTT}.
Suppose that $\Omega \subset \Rn$ is a bounded convex domain and that the evolving set is contained in
$\Omega$.
Furthermore, we need to assume that the one-homogeneous extension of $\beta$ as in
\eqref{one-homogeneous-extension} is convex.

The insight of Chambolle is to formulate the minimizing scheme of F.~Almgren, J.~E.~Taylor and L.~Wang
\cite{ATW} in terms of the signed distance function, so that the evolving set is its level set.
In \cite{Chambolle}, he proposed the time discretization by the minimization problem ($f \equiv 0$
  in
\cite{Chambolle})
\begin{align}
  \label{minimization}
  v_{m+1} \leftarrow \argmin_{v\in L^2(\Omega)} \pth{\frac{1}{2h} \norm{v - w_m}^2 + \int_\Omega \sigma(D v) \dx -
  \ang{f_m, v}},
\end{align}
where $h > 0$ is a chosen time step, $w_m$ is the signed distance function of the level set $\set{v_m < 0}$ at the previous time step $m$, induced by
the metric given by the mobility $\beta$ and $f_m = f(\cdot, mh)$. The minimization is performed over all $v \in
L^2(\Omega)$, and $\norm{\cdot}$ and $\ang{\cdot, \cdot}$ are the $L^2(\Omega)$-norm and inner product,
respectively.
The total variation energy $\int_\Omega \sigma(Dv) \dx$ is defined as the $L^2$ lower
semicontinuous envelope of the functional $v \mapsto \int_\Omega
\sigma(\nabla v) \dx$ defined for $v$ in the Sobolev space $W^{1,1}(\Omega)$. Note that the
functional in \eqref{minimization} is the Moreau-Yosida regularization
of the total variation energy with parameter $h$, and the minimization problem \eqref{minimization} is equivalent to
the resolvent problem for the total variation energy. In other words, it
is the implicit Euler discretization of the total variation flow \cite{ACM}. Since $\abs{\nabla v_{m+1}} \approx \abs{\nabla w_m} = 1$ when $h$ is small, we can deduce that $\frac{v_{m+1} - w_m}h \approx V$ using
\eqref{level-set-eq}.

The full minimizing movements algorithm to find a discrete sequence of approximations $E_m$, $m = 0, 1, 2, \ldots$ of the evolving set $\set{E_t}_{t \geq 0}$ at time steps $t_m
= m h$ reads: Set
$E_0$ as the initial data and then iteratively for $m = 0,
\ldots$ do
\begin{align}
  \label{Chambolle-algorithm}
  \begin{aligned}
  w_m &\leftarrow \signdist_\beta E_m\\
  v_{m+1} &\leftarrow \argmin_{v \in L^2(\Omega)} \pth{\frac{\mu}{2} \norm{v - w_m}^2 + \norm{\sigma(\nabla v)}_1 -
  \ang{f_m, v}},\\
  E_{m+1} &\leftarrow \set{v_{m+1} < 0},
  \end{aligned}
\end{align}
where $\mu = \frac 1h$, and $f$ is a given source, and we write $\norm{\sigma(\nabla v)}_1 =
\int_\Omega \sigma(Dv) \dx$. Note that the minimization is equivalent to
the minimization of
\[
  \frac{\mu}{2} \norm{v - (w_m + \mu^{-1} f_m)}^2 + \norm{\sigma(\nabla v)}_1.
\]

The signed distance function must correspond to the anisotropy $\beta$.
Recall that we assume that $\beta$ is one-homogeneously extended to $\Rn$ as in
\eqref{one-homogeneous-extension}, and such a extension is an anisotropy. In particular, $\beta$ is
assumed to be convex (but it needs not to be symmetric with respect to the origin).
We define the signed distance as
\begin{align*}
  \signdist_\beta E_m(x) &:=
  \inf_{y \in E_m} \beta^\circ(x - y) - \inf_{y \in E_m^\compl} \beta^\circ(y - x)\\
  &\ =
  \begin{cases}
    \inf_{y \in E_m} \beta^\circ(x - y), & x \notin E_m,\\
    - \inf_{y \in E_m^\compl} \beta^\circ(y - x), & x \in E_m,
  \end{cases}
\end{align*}
where $\beta^\circ(x) := \sup \set{x \cdot p: \beta(p) \leq 1}$ is the convex polar of $\beta$.
If $\beta(p) = |p|$, $\signdist_\beta$ is
just the standard signed distance function induced by the Euclidean metric.

Let us motivate the above choice of $w_m$.
If we set $w_m$ to be the signed distance above, we have $\beta(\nabla w_m) = 1$ a.e.,
see \cite{Rockafellar}.
Then
performing \eqref{Chambolle-algorithm}, $v_{m+1} - w_m$ is approximately $(\kappa - f_m) h$ and the free boundary
advances in the normal direction by
\begin{align*}
 -\frac{(\kappa - f_m) h}{|\nabla w_m|} = -\beta\pth{\frac{\nabla w_m}{|\nabla w_m|}}
 (\kappa - f_m) h = \beta(\nu) (-\kappa + f_m) h
\end{align*}
yielding the correct free boundary velocity.

As $h \to 0$, the evolution will converge to a continuous evolution
$\set{E_t}_{t\geq 0}$, see \cite{Chambolle}. Provided that there is no fattening, this evolution will be
the unique solution of the
anisotropic mean curvature flow \cite{CMP,CMNP,Ishii14} and, in the crystalline case in particular, the unique viscosity solution
solution of the crystalline mean curvature flow \cite{GP16,GP_CPAM}.

It might seem that the minimization problem in \eqref{Chambolle-algorithm} is rather difficult for
numerical computation, mainly due to the non-differentiable second term. In particular, the standard
minimization methods like conjugate gradients or Newton iteration are poorly suited.
For this reason, Chambolle proposed an iterative algorithm in \cite{Chambolle}.
More recently, it was recognized in \cite{OOTT} that the minimization problem can be addressed
by the so-called proximal algorithms \cite{ParikhBoyd}, the alternative direction method of multipliers (ADMM) or the split Bregman
method \cite{GO}.
 To find the minimizer $v$ of $\frac{\mu}{2} \norm{v - u}^2 + \norm{\sigma(\nabla
v)}_1$, we choose $\lambda > 0$, set $b_0 = d_0 = 0$ and then iterate for $k = 0, 1, \ldots$
\begin{subequations}
  \label{split-bregman}
\begin{align}
  \label{minimization-v}
  v_{k+1} &\leftarrow \argmin_{v} \frac{\mu}{2} \norm{v - u}^2 + \frac{\lambda}{2} \norm{d_k -
  \nabla v - b_k}^2,\\
  \label{minimization-d}
  d_{k+1} &\leftarrow \argmin_{d} \norm{\sigma(d)}_1 + \frac{\lambda}{2} \norm{d - \nabla v_{k+1} -
b_k}^2,\\
\label{bregman-iteration}
  b_{k+1} &\leftarrow b_k + \nabla v_{k+1} - d_{k+1},
\end{align}
\end{subequations}
until some stopping condition is reached, typically when $\norm{v_{k+1} - v_k}_2$ is
sufficiently small.
Heuristically, this scheme introduces a new gradient variable $d$, and then enforces
the constraint $d = \nabla v$ by a quadratic penalty.
Since we are minimizing a sum of convex terms over two variables, the problems can be decoupled
into iterating \eqref{minimization-v} and \eqref{minimization-d}.
\eqref{bregman-iteration} is called a Bregman iteration, and it
helps to enforce the constraint \emph{exactly}.
Note that when convergence is achieved,  \eqref{bregman-iteration} implies $d = \nabla v$. For a
detailed discussion of the motivation, convergence and other properties, see \cite{GO}.

The advantage of this iteration process is the simplicity of the subproblems. The first
minimization problem \eqref{minimization-v} is equivalent to finding the solution $v$ of
\begin{align}
  \label{elliptic-problem}
  (\mu - \lambda \Delta) v = \mu u + \lambda \divo (b_k - d_k) \qquad \text{in $\Omega$},
\end{align}
with an appropriate boundary condition, for instance Neumann.
It is not necessary to solve it accurately, so one or two Gauss-Seidel iterations are sufficient
\cite{GO}.

To find a numerical solution, we use the finite difference method (FDM)  or the finite element method
(FEM) to discretize
\eqref{Chambolle-algorithm}, see Section~\ref{sec:discretization}. In both cases, we represent $d,
b$ by discrete values on the lattice
(FDM) or on the elements (FEM).
This has the important consequence that the minimization problem \eqref{minimization-d} for $d$ in
the discrete case completely decouples. Then, for each node or element $i$
the $i$th component of the minimizer $d_{k+1,i}$ is given by the so-called shrink operator \cite{GO}
\begin{align}
  \label{d-as-shrink}
  d_{k+1,i} = \shrink_\sigma((\nabla v_{k+1} + b_k)_i, 1/ \lambda).
\end{align}
Note that the shrink operator can be expressed using the orthogonal projection on the Wulff shape
$\Wulff$ of $\sigma$ \cite{OOTT},
\begin{align*}
  \shrink_\sigma(\xi, 1 / \lambda) := (I - P_{\Wulff/\lambda}) (\xi).
\end{align*}
In typical cases of isotropic, cubic and hexagonal anisotropies, that is, when the Wulff shape $\Wulff$ is a
sphere, a cube or a hexagonal prism, respectively, the orthogonal projection is very simple. More general Wulff
shapes can be handled by the method proposed in \cite{OOTT}.

It is interesting to relate the Bregman iterate $b$ to the Cahn-Hoffman vector field in the
definition of the anisotropic (crystalline) mean curvature.
In particular, for a node or an element $i$, if convergence is achieved, \eqref{d-as-shrink} yields $d_i =
\shrink_\sigma(d_i + b_i, 1 / \lambda)$ and we have
  $b_i = P_{\mathcal W/\lambda}(d_i + b_i)$.
From the last equality, we see that either $d_i = 0$, but then $\lambda b_i \in \Wulff$, or
$d_i
\neq 0$, but then $\lambda b_i \in \partial \Wulff$ and $d_i$ is a normal of $\partial
\Wulff$ at $\lambda b_i$. In other words, $\lambda b_i \in \partial \sigma(d_i)$. From
\eqref{elliptic-problem} we deduce that $\frac{v- u}h = \lambda \divo b$.
Hence $\lambda b$ is a discrete Cahn-Hoffman vector field for the resolvent problem
\eqref{minimization}.

\subsection{Discretization: FDM vs FEM}
\label{sec:discretization}

\begin{figure}
  \begin{minipage}[c]{.5\textwidth}
  \begin{subfigure}{\textwidth}
    \includegraphics[width=\textwidth]{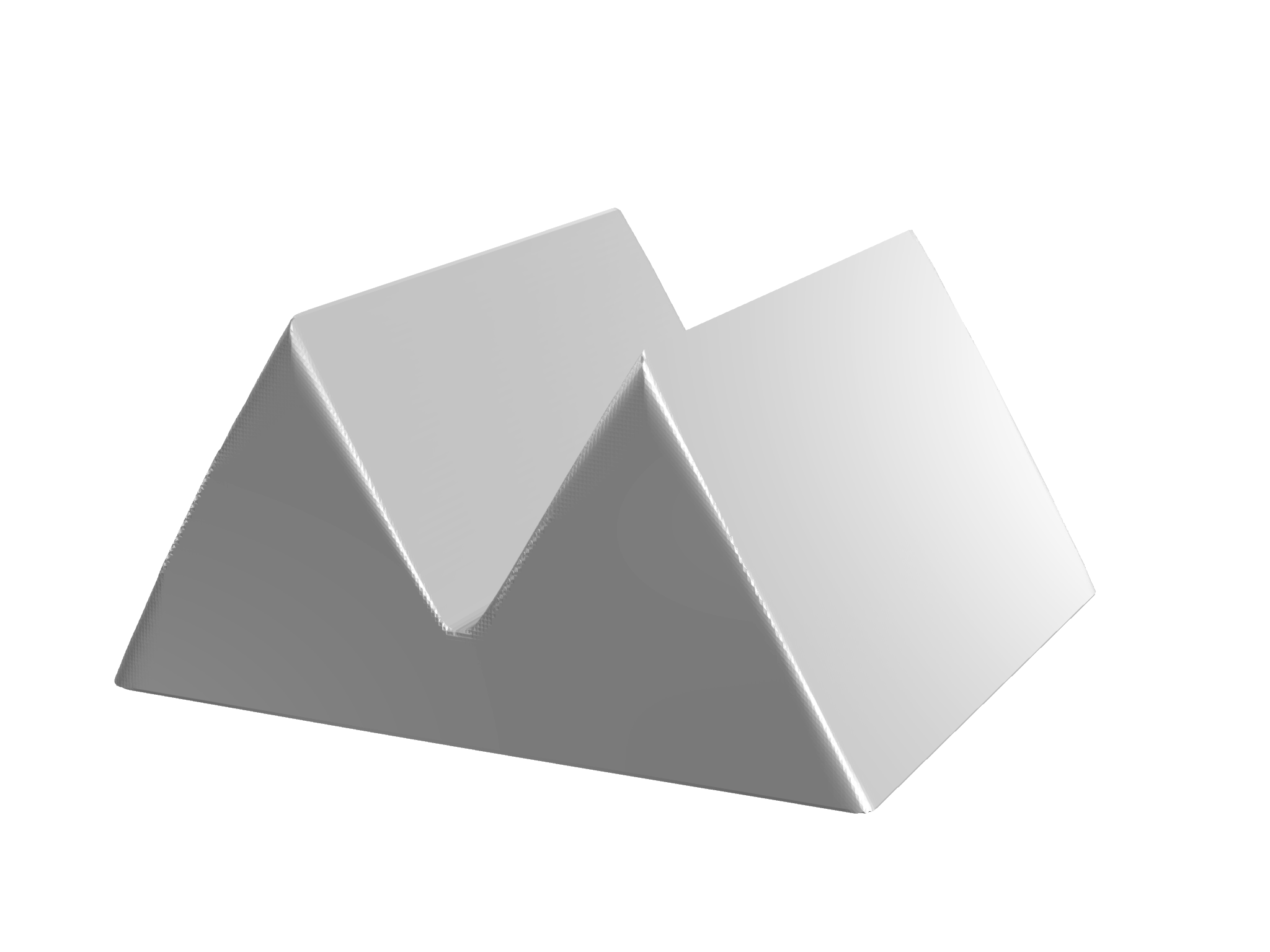}
  \end{subfigure}
  \end{minipage}
  \begin{minipage}[c]{.33\textwidth}
    \begin{subfigure}{0.48\textwidth}
    \includegraphics[width=\textwidth]{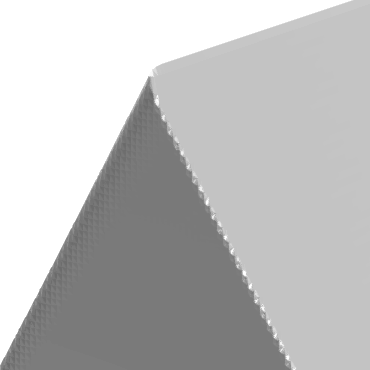}
  \end{subfigure}
  \begin{subfigure}{0.48\textwidth}
    \includegraphics[width=\textwidth]{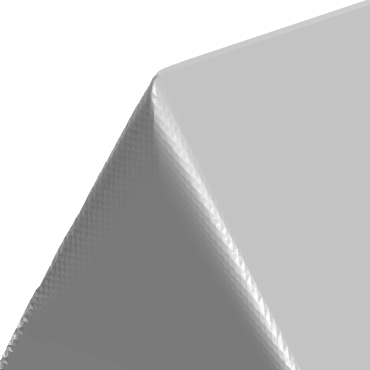}
  \end{subfigure}\\[0.2cm]
  \begin{subfigure}{0.48\textwidth}
    \includegraphics[width=\textwidth]{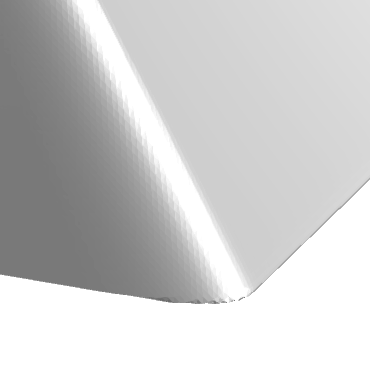}
  \end{subfigure}
  \begin{subfigure}{0.48\textwidth}
    \includegraphics[width=\textwidth]{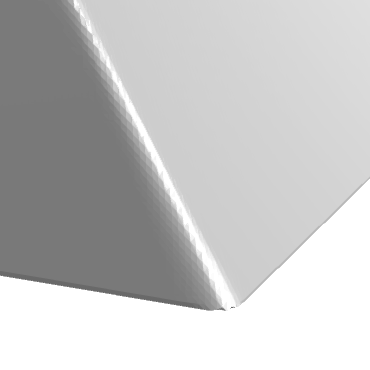}
  \end{subfigure}
  \end{minipage}
  \caption{Rounding of edges in the computation in Figure~\ref{fig:two_tris}; neighborhood of the
    top left vertex (top row) and bottom right vertex (bottom row) are magnified: FDM (left
  column) vs. FEM (right column).}
  \label{fig:rounding_threed}
\end{figure}

The standard way \cite{GO,OOTT} to discretize \eqref{split-bregman} is to use the finite difference
method (FDM). This is quite straightforward, see also \cite{Chambolle}, and  seems sufficient in
two dimensions and for applications in image processing.
For instance, \eqref{elliptic-problem} is discretized using the standard central difference scheme
on a $2n+1$ point stencil.
However, in three dimensions this approach introduces unwanted artifacts, for instance rounding of edges in
certain directions, see Figure~\ref{fig:rounding_threed}.
This is most likely caused by the fact that the gradient $d$ and the Cahn-Hoffman vector field $b$
on a given cubic cell are given only by four out of the eight nodes. This restriction of the number of
degrees of freedom for the gradient limits the ability of the discretization to account for a rapidly changing vector
field near an edge.

We propose a finite element method (FEM) discretization. We minimize the functional in
\eqref{minimization-v} on the space of piece-wise linear functions ($P^1$~elements) on a tetrahedral mesh.
The vector fields $\nabla v$, $d$ and $b$
are then approximated by piece-wise constant vector fields ($P^0$~elements) on this mesh. The minimization problem
\eqref{minimization-d} on the space of piece-wise constant ($P^0$~elements) completely decouples on each element and the minimizer is given by \eqref{d-as-shrink}.

\begin{figure}
  \includegraphics[trim={1in 1.1in 1in 0.4in},clip,width=\textwidth]{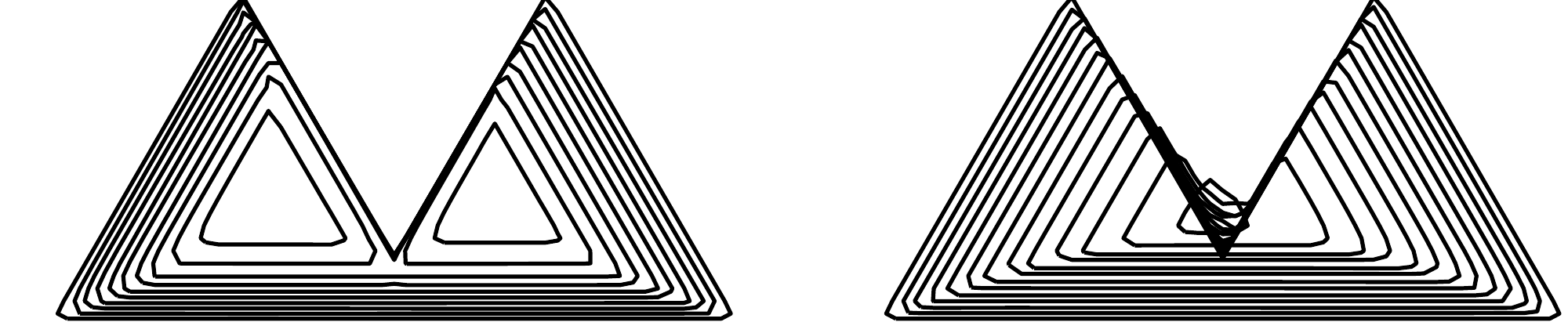}
  \caption{A triangular anisotropy demonstrates a topological change in 2D for non-even
    anisotropies. Left: FDM Right: FEM. Parameters $h = 10^{-4}$, $M = 64$, and plot
    step $0.005 = 50h$.}
  \label{fig:triangles}
\end{figure}

The merit of this approach is the six-fold increase of the number of degrees of freedom of $d$ and $b$
(for the cost of also increasing the memory and computational time requirements),
which allows for a better resolution of edges, see Figure~\ref{fig:rounding_threed}.
Additionally, it is consistent with our redistance approach in Section~\ref{sec:redistance}.
However, at some extreme cases in two dimensions the FEM performs worse than FDM, see
Figure~\ref{fig:triangles} for an example.

\begin{figure}
  \centering
  \fig{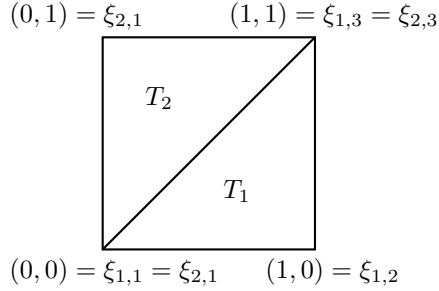}{4in}
  \caption{Tessellation of the unit cube in $n = 2$.}
  \label{fig:mesh}
\end{figure}

Let us present the mesh construction in an arbitrary dimension to fix the idea.
We assume that the computational domain is a cube $\Omega$ is subdivided into $M^n$ smaller
cubes of equal size, $M \in \N$ is the resolution.
We split each of the smaller cubes into $n!$ simplices in the following way. By translation and
scaling, we may assume that the small cube is $Q = [0,1]^n$.
Let $\set{p_m}_{m=1}^{n!}$ be the $n!$ permutations of $(1, \ldots, n)$. For each $m = 1, \ldots,
n!$, we define
the sequence of $n + 1$ vertices $\xi_{m, 1}, \ldots, \xi_{m, n+1}$ of $Q$ using
\begin{align*}
  \xi_{m, 1} &= (0, \ldots, 0),\\
  (\xi_{m, j+1})_i &=
\begin{cases}
  (\xi_{m,j})_i + 1 & \text{if $i = (p_m)_j$,}\\
  (\xi_{m,j})_i & \text{otherwise},
\end{cases}
\qquad i = 1, \ldots, n.
\end{align*}
Then their convex hulls,
\begin{align*}
  \operatorname{conv}(\xi_{m,1}, \ldots, \xi_{m, n+1}), \qquad 1 \leq m \leq n!,
\end{align*}
are simplices that form a tessellation of $Q$, see Figure~\ref{fig:mesh}.

We use the mesh that is given by the above tessellation of each of the $M^n$ cubes in $\Omega$. The
total number of elements is $n! M^n$.

\subsection{Redistance}
\label{sec:redistance}
At each time step of the algorithm \eqref{Chambolle-algorithm}, the signed distance function $\signdist_\beta \set{v_m < 0}$ to the
$0$-level set of the
solution of the minimization problem has to be recomputed. (This
is not strictly true, see the discussion in Section~\ref{sec:implementation-notes}.) It is sometimes referred to as
\emph{redistance}. The problem amounts to solving the (anisotropic) eikonal equation $|\nabla
w|_\beta := \beta^\circ(\nabla w) = 1$ with boundary data $w = 0$ on $\partial E$. There are various efficient methods for doing
this, including an iteration scheme \cite{SSO}, as well as more direct algorithms like the fast marching
method \cite{Sethian96} or the fast sweeping method \cite{Zhao}. We choose the fast sweeping method
due to its simplicity and efficiency.

The boundary $\partial E$ is given as the $0$-level set of a function $v$ with discrete values on
a regular grid.
Unfortunately, in general the $0$-level set of the resulting approximation $w$ of the signed
distance function will be different from the $0$-level set of the original function $v$. The fast
marching and fast sweeping methods require an initialization step, where the distance function is
assigned at grid points that are direct neighbors of the $0$-level set of $v$. A special care must
be taken so that the interface is not moved unnecessarily, as these effects might quickly
accumulate over a series of consecutive time steps. This is especially important at points where the surface
should not move, a typical case for non-convex and non-concave facets, where this unwanted
redistance effect might dominate
the evolution.

There are a few standard schemes to initialize the nearby values in the literature, by
analyzing the intersection of the level set with the grid lines \cite{AS,Chambolle}. However, they
do not seem to reproduce the correct value of the distance function even if the level set is flat,
which is a common situation in the crystalline mean curvature flow.
Furthermore, the generalization to three dimensions seems unnecessarily complicated.

We choose a naive method that appears to be superior in our case, is very simple to implement in
an arbitrary dimension, and that computes the \emph{exact} signed distance function in a
neighborhood of the flat facets (away from vertices and edges). The idea is to split the $M^n$ cubes of the uniform grid
into $n!$ simplices as explained in the construction of the mesh for the
finite element method in Section~\ref{sec:discretization} and suppose that $v$ is affine on each of
them. This is very much in the spirit of the marching tetrahedra method \cite{DK}. Then for all elements
that the $0$-level set intersects, that is, on which $v$ changes sign, we set the initial
value of the signed distance function at each vertex of the element to be the value of $v$
normalized by the $\beta$-norm of the gradient of the affine function given by $v$ on the element. If a given grid node
is a vertex of multiple elements intersecting the $0$-level set, we set the initial value to be the
minimum over all the elements. To be more explicit, suppose that $v_i$, $w_i$, $i = 1, \ldots, N$
are the values at the grid nodes $x_i$, and $\mathcal T_j$, $j = 1, \ldots, K$ are the simplices on which
$v$ changes sign. We initialize $w_i$ to
\begin{align*}
  w_i = (\sign{v_i}) \inf_{\substack{1 \leq j \leq K\\x_i \in \partial \mathcal T_j}}
  \frac{|v_i|}{\beta^\circ(\at{\nabla v}{\mathcal
T_j})},
\end{align*}
where the infimum is defined as $+\infty$ if it is over an empty set and $\at{\nabla v}{\mathcal
T_j}$ is the (constant) value of $\nabla v$ on the simplex $\mathcal T_j$. This initialization method is
second order accurate, $O(M^{-2})$, near a smooth surface, in contrast to the first order accuracy of the
initialization in \cite{AS,Chambolle}. Moreover, the unwanted artifacts caused by the redistance seem
to be reduced, for instance compare the preservation of non-convex/non-concave facets in
Figure~\ref{fig:many_squares}, even for a modest space resolution $M=64$.

After this initialization step, we perform the $2^n$ sweeps of the fast sweeping method \cite{Zhao}.

\subsection{Notes on the implementation}
\label{sec:implementation-notes}

Let us give a few notes on our implementation from a practical point of view.

In the discussion of the iteration \eqref{split-bregman}, we suggested to initialize $b_0 = d_0 =
0$. However, this is unnecessary since the iteration converges to the unique minimizer no matter
the initial guess. Since the vector fields $b$ and $d$ vary relatively slowly from one time step to
another, we can reuse the value of $b$ and $d$ from the previous time step to start the iteration. The
main consequence is that the number of iterations of \eqref{split-bregman} necessary to obtain a
reasonably accurate result is dramatically decreased as the work is spread out over consecutive time steps, see Table~\ref{tab:performance}.

The choice of the stopping criterion for the iteration \eqref{split-bregman} is important. In our
implementation we stop once $\norm{v_{k+1} - v_k} < \e_{\rm btol}$ for some given $\e_{\rm btol} > 0$, where
$\norm{\cdot}$ is for example the $\ell^2$ norm.
When $\e_{\rm btol}$ is chosen too large, the facets do not become flat as there are too few iterations
to reduce the lowest frequency mode
of the error, whose wavelength is proportional to the size of the facets. On the other hand, if $\e_{\rm btol}$ is
chosen too small, the number of iterations becomes unnecessarily high for no gain in accuracy, which is limited by
the time step error $O(h)$ and the ability to resolve the corners and edges of the crystal with
accuracy $O(M^{-1})$. Unfortunately, there does not seem to be an explicit way to estimate the
necessary $\e_{\rm btol}$ and some experimentation is required, see Figure~\ref{fig:lambda-analysis}.

Another source of numerical problems is the computation of the distance function in
\eqref{Chambolle-algorithm}, see the discussion in Section~\ref{sec:redistance}.
However, often $v_m$ is a very good approximation of $w_m$ close to the zero-level set and hence we can take $w_m = v_m$ and
skip the distance computation for a few time steps, after which we need to recompute the distance
function. This can help reducing the redistance artifacts that sometimes appear, especially if the time step is very
small and there are parts of the surface that do not move, such as non-convex/non-concave facets or
curved sections. We did not use this in the results presented here.

\subsection{Numerical results}
\label{sec:numerical-results}

To illustrate the performance, we present a few simple numerical results based on our implementation of the above algorithm in the Rust programming language. The domain is always taken to be $\Omega = (-\frac12, \frac12)^n$ and $\lambda = \frac\mu8$ in \eqref{minimization-v}, \eqref{minimization-d}, see Section~\ref{sec:choice-of-lambda}.  The stopping condition is chosen as $\norm{v_{k+1} - v_k}_{\ell^2} < 10^{-5}$ for $n = 2$ and $\norm{v_{k+1} - v_k}_{\ell^2} < 10^{-4} M^{\frac12}$ for $n = 3$, where $\norm{\cdot}_{\ell^2}$ is the discrete $\ell^2$-norm, so that the average stopping tolerance per element scales like $M^{-1}$. This choice seems to give reasonable results across all computations presented here. We use $M = 64$ in $n = 2$ and $M = 256$ for $n = 3$. The performance for $n=3$ is $\approx 1{\rm min}/\rm{timestep}$ for FDM and $\approx 5{\rm min}/\rm{timestep}$ for FEM on a single core of Intel Core i7-4770K @ 3.50 GHz. The figures depicting the solution show the zero-level set of the numerical solution $v_m$.

\subsubsection{Estimating the error of solutions}
To quantify the error between the numerical and the exact solution, we consider the Hausdorff distance between two surfaces $\Gamma_1, \Gamma_2 \subset \Rn$ in the maximum and $L^2$ norms as
\begin{align*}
\dist_{H,\infty}(\Gamma_1, \Gamma_2) &:= \max\pth{\max_{x \in \Gamma_1} \dist(x, \Gamma_2),\max_{x \in \Gamma_2} \dist(x, \Gamma_1)},\\
\dist_{H,2}(\Gamma_1, \Gamma_2) &:= \pth{\int_{\Gamma_1} \dist(x, \Gamma_2)^2 \diff{\mathcal H^{n-1}} + \int_{\Gamma_2} \dist(x, \Gamma_1)^2 \diff{\mathcal H^{n-1}}}^{\frac12}.\\
\end{align*}
Note that $\dist_{H,\infty}$ is the usual Hausdorff distance.
Both are defined to be $\infty$ if exactly one of the surfaces is empty.

\begin{figure}
  \includegraphics[width=\textwidth]{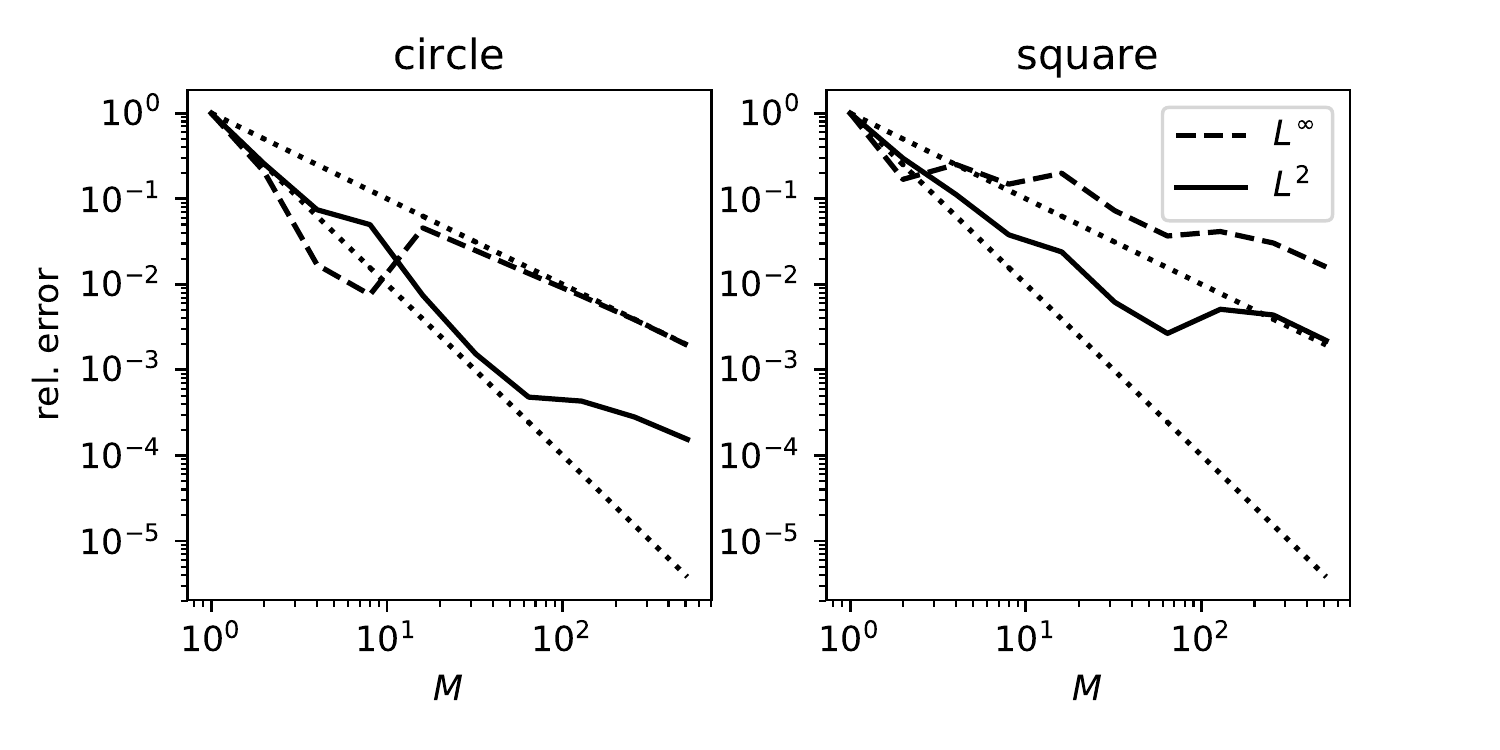}
  \caption{Relative error of the numerical Hausdorff distance for two concentric circles of radii $0.3$ and $0.4$ (left) and two concentric squares of side-length $0.6$ and $0.8$ (right) as a function of the mesh resolution $M$. The distance function in our implementation of the fast sweeping method is only first-order accurate, which influences the computation for larger $M$. This is more apparent in the case of the maximum distance. Dotted lines are $M^{-1}$ and $M^{-2}$ for comparison.}
  \label{fig:distance_error}
\end{figure}

Numerically, these two distances between two level sets are computed with the discrete distance function given by the fast sweeping method with the initialization described in Section~\ref{sec:redistance}. The integrals and maxima are approximated by those of the piece-wise linear distance function over a mesh generated from the level set function by the marching tetrahedra method \cite{DK} in two or three dimensions.
Both distances work well for smooth surfaces, but $\dist_{H,\infty}$ penalizes corners and overemphasizes errors on a small part of the boundary. In Figure~\ref{fig:distance_error} we test the accuracy of this numerical algorithm in cases of the level set of $(x_1^2 + x_2^2)^{1/2} - r$ (circle) and $\max(|x_1|, |x_2|) - r$ (square) with $r = 0.3, 0.4$ for each. Initially, when $M^{-1} \approx 0.1$, the accuracy is dominated by the accuracy of the level set reconstruction using the marching tetrahedra and the distance function initialization explained in Section~\ref{sec:redistance} and both of these are second-order accurate. Once $M^{-1} \ll 0.1$, the first-order discretization error in the fast sweeping algorithm dominates. Therefore we report the errors using $\dist_{H,2}$, which appears to be less sensitive to this error.

\subsubsection{Choice of $\lambda$}\label{sec:choice-of-lambda}
\begin{figure}
  \includegraphics[width=\textwidth]{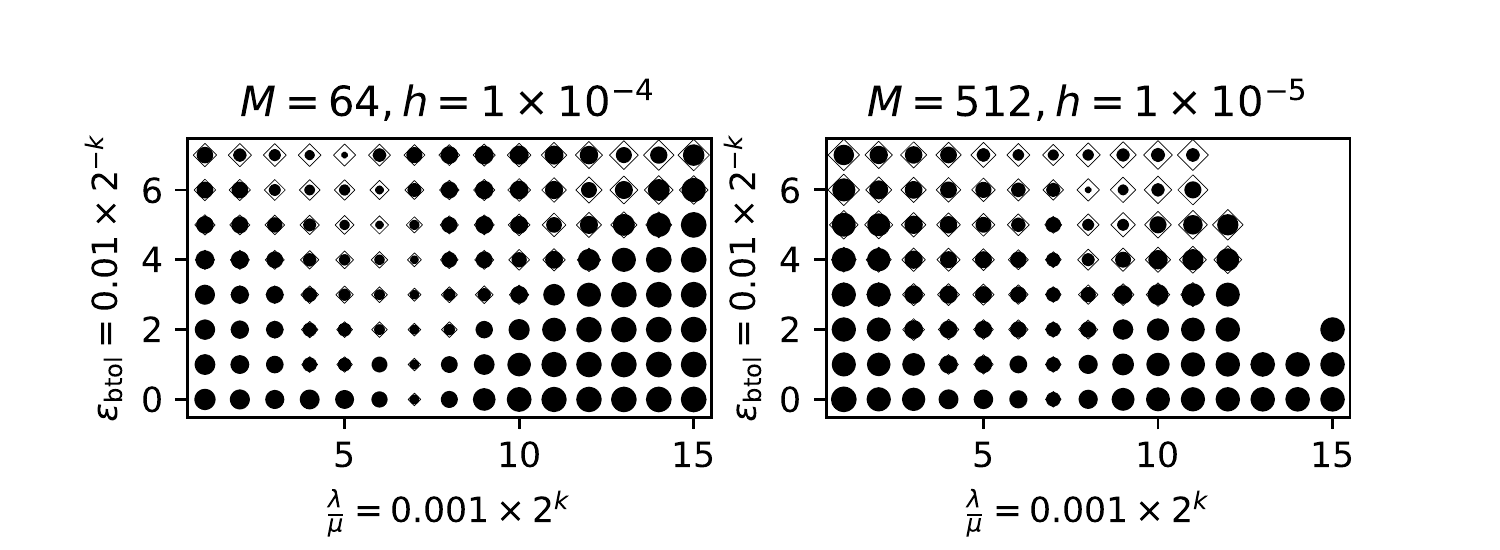}
  \caption{Performance for different $\frac \lambda\mu$ and $\e_{\rm btol}$ for the shrinking hexagon in $n = 2$. The area of the discs is proportional to the logarithm of the error, while the area of the diamonds is proportional to the logarithm of the computational time. $\frac{\lambda}{\mu} = 0.001 \times 2^7 \approx \frac 18$ is generally well-performing.}
  \label{fig:lambda-analysis}
\end{figure}

Let us address the choice of the parameter $\lambda$ in proportion to $\mu$. Theoretically, the split Bregman iteration will converge for any $\lambda > 0$. However, this parameter significantly affects the speed of convergence. In the original paper \cite{GO} it was suggested to use $\lambda = 2 \mu$. But after somewhat extensive testing in our case it appears that the value $\lambda = \frac{\mu}{8}$ gives the fastest convergence and smallest error overall over both two and three dimensional computations. In this test, we simply ran the computation for many choices of parameters $\frac\lambda\mu$ and stopping conditions $\e_{\rm btol}$ in $\norm{v_{k+1} - v_k}_{\ell^2} < \e_{\rm btol}$ for a self-similar shrinking hexagon in two dimensions and shrinking cube in three dimensions with various $M$ and $h$, and plotted the errors and computational time. See Figure~\ref{fig:lambda-analysis} for an example. Due to the computational cost in three dimensions, the testing was much less exhaustive there.

\subsubsection{Self-similar solutions}
\label{sec:self-similar-tests}
\begin{figure}[p]
  \includegraphics[width=\textwidth]{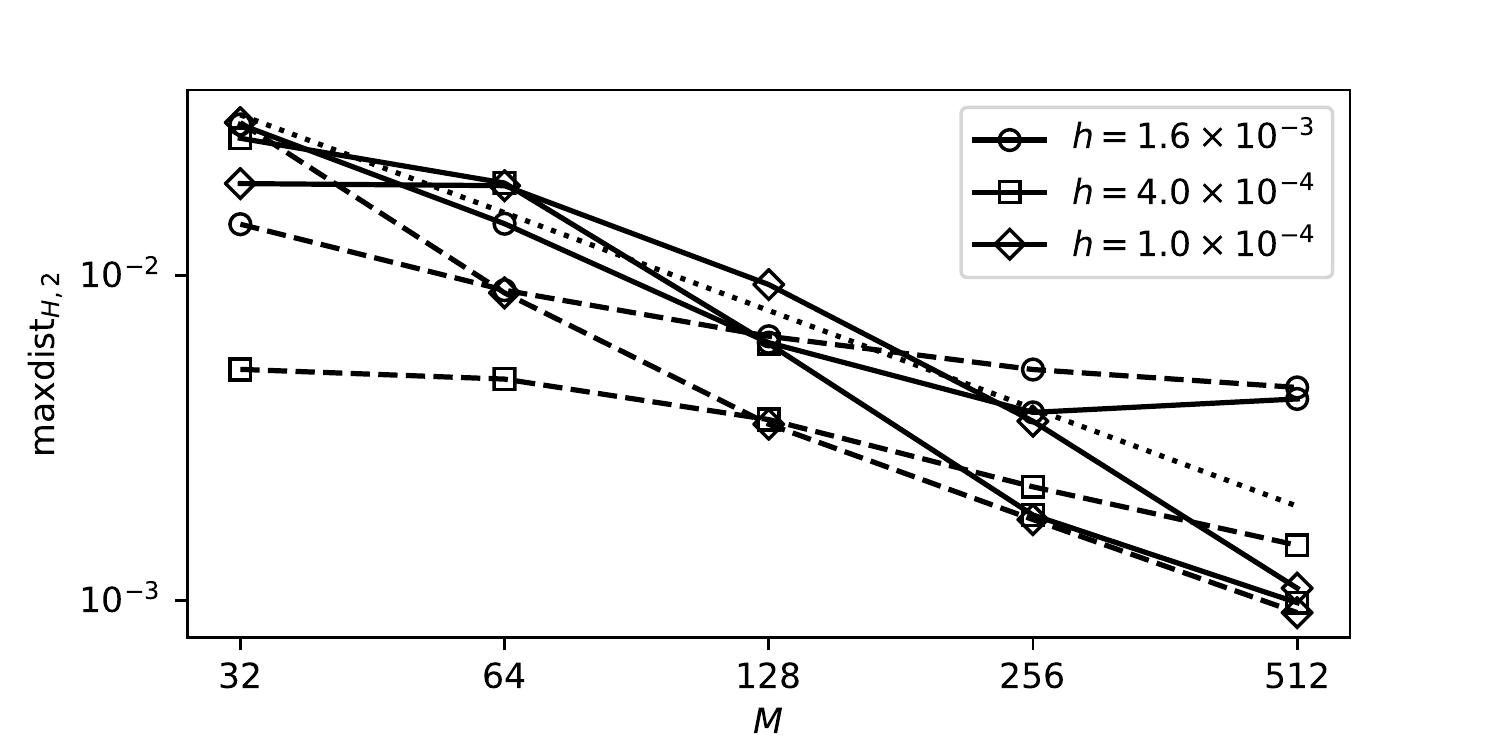}
  \caption{The error of a self-similar shrinking hexagonal Wulff shape of side-length $0.4$ in $n=2$ with the numerical solution using the FEM (solid) and FDM (dashed) discretization (maximum over $t \in [0, 0.05]$). Dotted line is $M^{-1}$.}
  \label{fig:hexagon_errors}
\end{figure}

\begin{figure}[p]
  \includegraphics[width=0.92\textwidth]{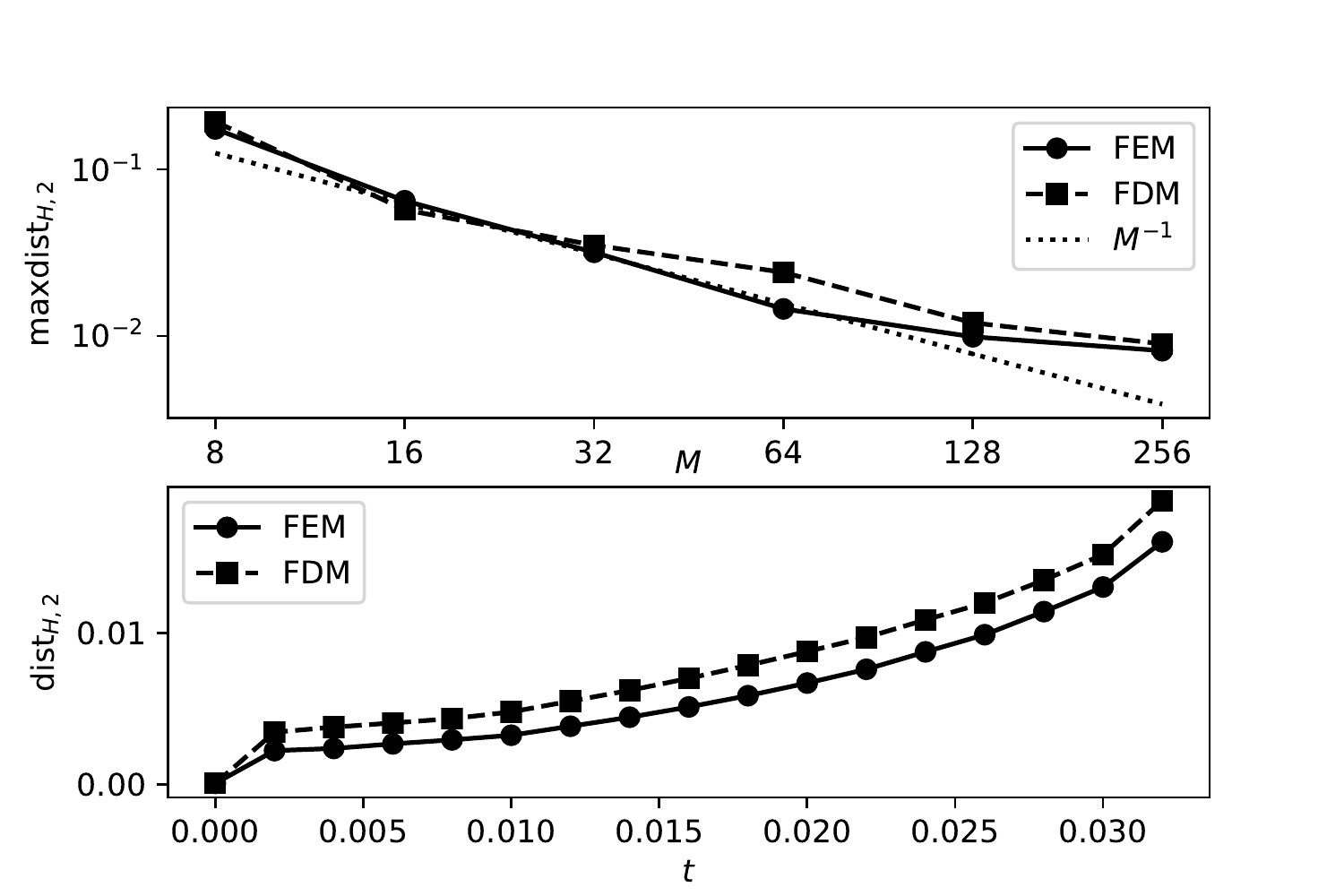}
  \caption{The error of a self-similar shrinking hexagonal Wulff shape with a hexagonal base with edges of length $0.4$ and height $0.8$ in $n=3$ with $h = 10^{-4}$ as the function of $M$ (top, maximum over $t \in [0, 0.026]$) and $t$ for $M = 128$ (bottom).}
  \label{fig:hexagonal_prism_errors}
\end{figure}

\begin{figure}[p]
  \includegraphics[width=0.92\textwidth]{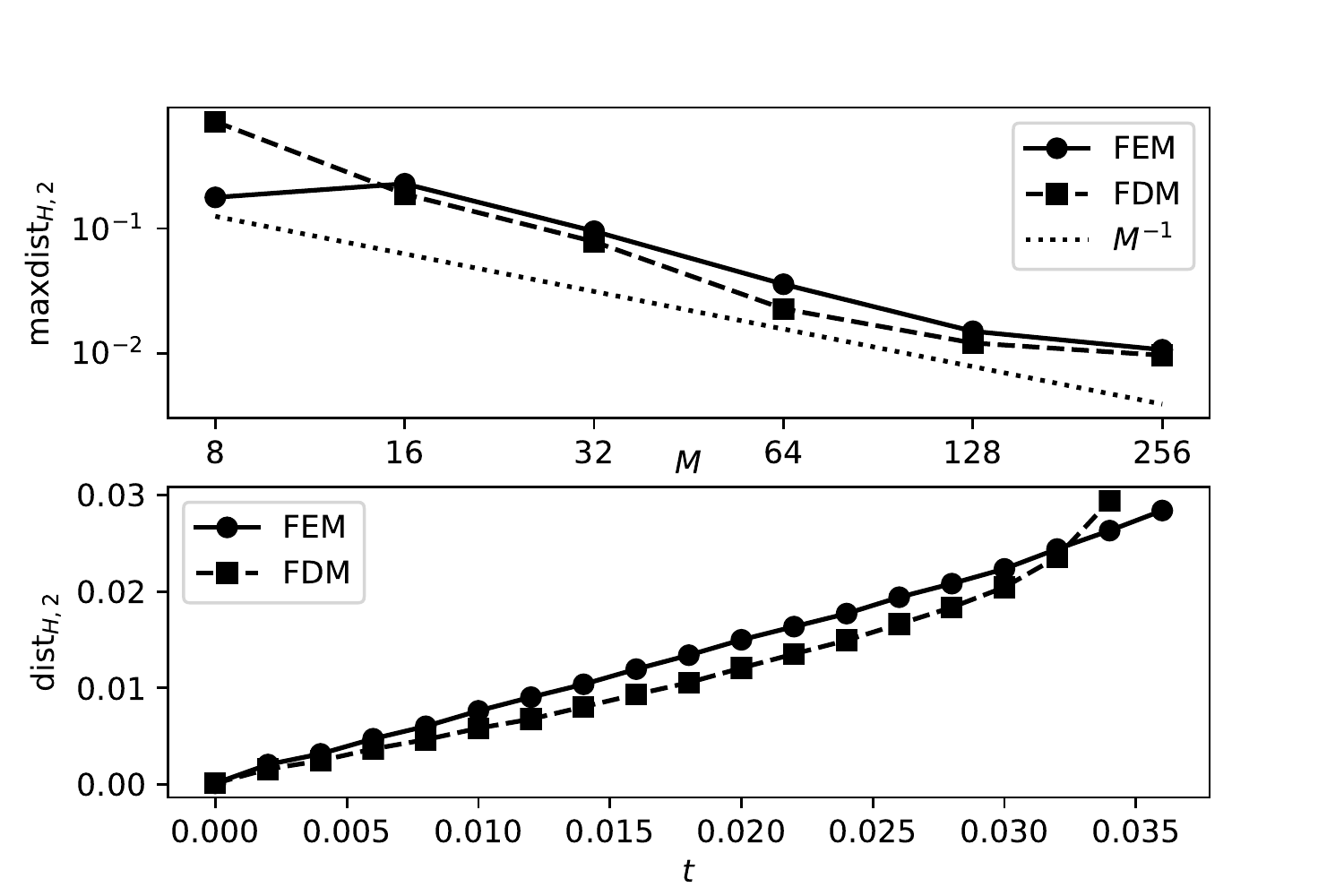}
  \caption{The error of a self-similar crystalline shrinking doughnut with $\ell^1$ anisotropy from Section~\ref{sec:shrinking-doughnut} of side-length $0.8$ in $n=3$ with $h = 10^{-4}$ as the function of $M$ (top, maximum over $t \in [0, 0.02]$) and $t$ for $M = 128$ (bottom).}
  \label{fig:torus_errors}
\end{figure}

\begin{figure}[p]
  \includegraphics[width=0.92\textwidth]{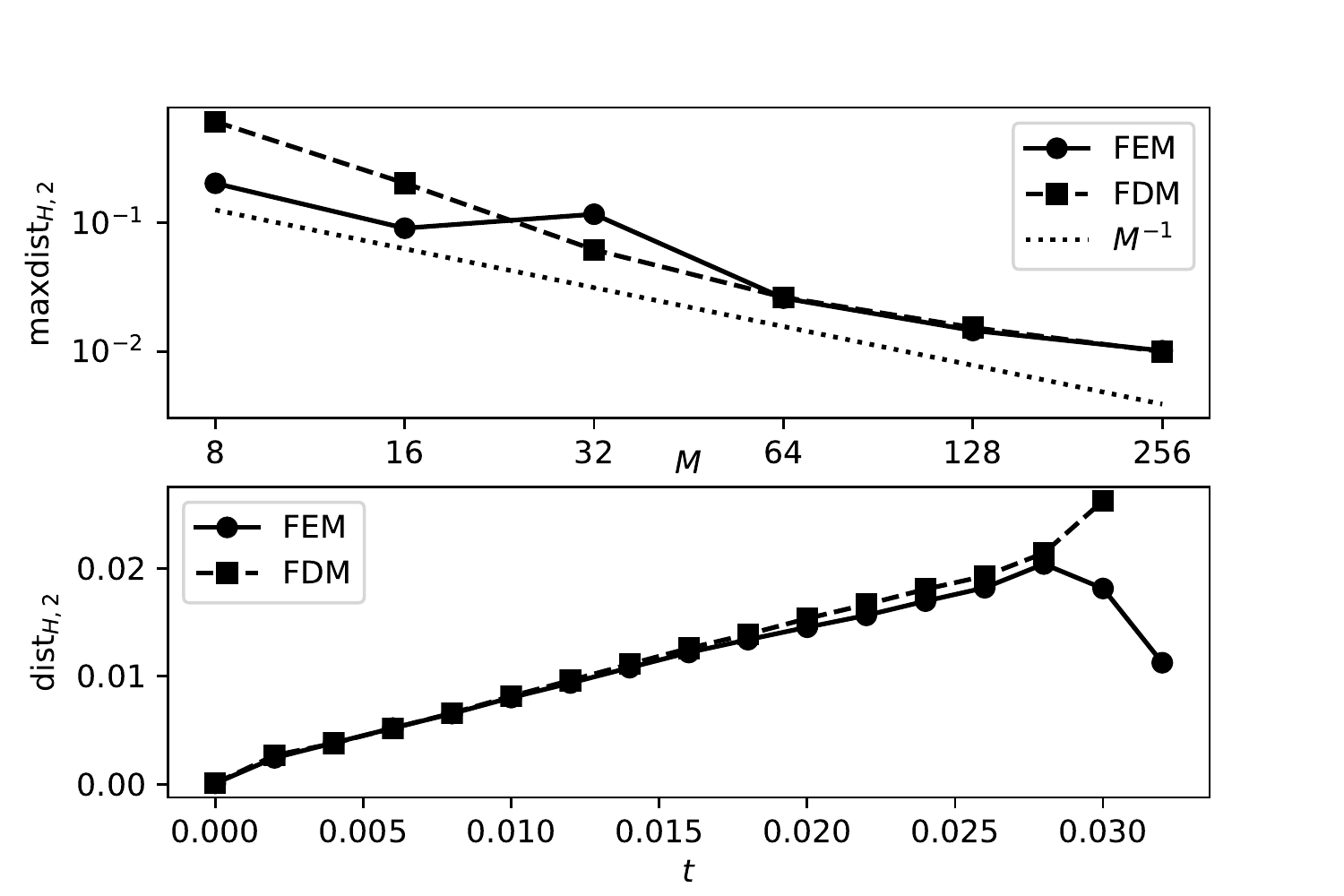}
  \caption{The error of a self-similar crystalline shrinking doughnut with hexagonal anisotropy from Section~\ref{sec:shrinking-doughnut} of hexagonal base with side-length $0.4$ and height $0.8$ in $n=3$ with $h = 10^{-4}$ as the function of $M$ (top, maximum over $t \in [0, 0.02]$) and $t$ for $M = 128$ (bottom).}
  \label{fig:hex_torus_errors}
\end{figure}

\begin{figure}[p]
  \includegraphics[width=0.92\textwidth]{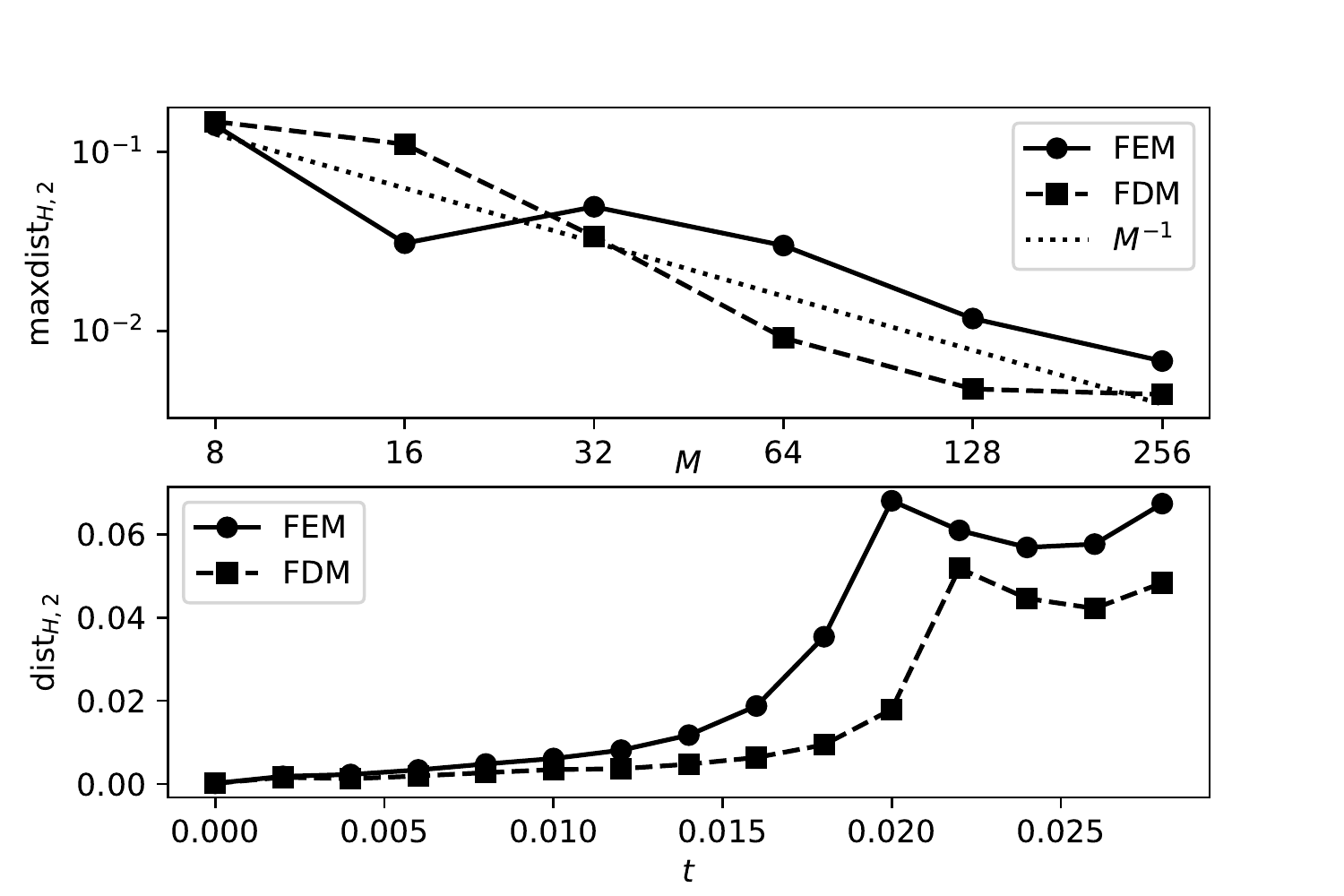}
  \caption{The error of a self-similar shrinking sponge from Section~\ref{sec:sponge} with $R_0 = 0.4$ in $n=3$ with $h = 10^{-4}$ as the function of $M$ (top, maximum over $t \in [0, 0.014]$) and $t$ for $M = 128$ (bottom).}
  \label{fig:sponge_errors}
\end{figure}

We test the method in two and three dimensions using self-similar solutions of $V = \beta(\nu) \kappa_\sigma$. In two dimensions, we use the shrinking Wulff shape solution with hexagonal symmetry: the Wulff shape is the regular hexagon with edge of length $1$, Figure~\ref{fig:hexagon_errors}. The results are given for three different time step sizes $h$. For implementation reasons, this and the below tests with hexagonal mobility have time scaled by the factor $\frac{2\sqrt{3}}3$. Therefore the extinction time in $n=2$ is not $t^* = \frac{0.4^2}2 = 0.08$ as expected but $\frac{0.4^2\sqrt3}4 \approx 0.069$, and $\frac{0.4^2\sqrt3}8 \approx 0.035$ in $n=3$. In three dimensions, we use the shrinking Wulff shape solution when the Wulff shape is the hexagonal prism with height $2$ and regular hexagonal base with edges of length $1$, Figure~\ref{fig:hexagonal_prism_errors}. Finally, we use the self-similar solutions constructed in Section~\ref{sec:self-similar}. Due to their instability, presence of different types of facets and simplicity, they seem to provide a convenient benchmark. For the crystalline shrinking doughnut, we use cubic and hexagonal anisotropies $\tilde\sigma$, see Figure~\ref{fig:torus_errors} and Figure~\ref{fig:hex_torus_errors} respectively. The error of the sponge solution is presented in Figure~\ref{fig:sponge_errors}. In all of these, the error is estimated as $\max_{m \in \mathcal M} \dist_{H,2} (\partial E(t_m), \partial E_m)$ where $\mathcal M$ is a set of selected time steps such that for all $m \in \mathcal M$ we have $t_m = 0.002z$ for some $z \in \Z$, and that $t_m$ are far enough from the extinction time $t^*$, where all the solutions diverge significantly due to the instability. Table~\ref{tab:performance} shows the details on the performance of the algorithm for $n = 3$ and $M = 64, 256$.

\begin{table}
\centering
\begin{tabular}{lrrrr}
\hline
\multicolumn{5}{c}{$M = 64$}\\
\hline
& \multicolumn{2}{c}{FDM} & \multicolumn{2}{c}{FEM}\\
\hline
Test & \#Breg/step & time/Breg (s) & \#Breg/step & time/Breg (s)\\
\hline
Hex Wulff & $24$ & $0.017$ & $29$ & $0.066$\\
$\ell^1$ doughnut & $35$ & $0.017$ & $52$ & $0.068$\\
Hex doughnut & $38$ & $0.017$ & $39$ & $0.075$\\
Sponge & $28$ & $0.014$ & $33$ & $0.044$\\
\hline
\hline
\multicolumn{5}{c}{$M = 256$}\\
\hline
& \multicolumn{2}{c}{FDM} & \multicolumn{2}{c}{FEM}\\
\hline
Test & \#Breg/step & time/Breg (s) & \#Breg/step & time/Breg (s)\\
\hline
Hex Wulff & $42$ & $0.93$ & $42$ & $2.8$\\
$\ell^1$ doughnut & $100$ & $0.68$ & $109$ & $3.0$\\
Hex doughnut & $93$ & $0.78$ & $98$ & $3.5$\\
Sponge & $133$ & $0.69$ & $125$ & $3.1$\\
\hline
\end{tabular}
\caption{Performance of the numerical method in $n=3$ with $M = 64$ and $M = 256$ for tests in Section~\ref{sec:self-similar-tests} based on the average number of Bregman iteration per time step (\#Breg/step) and average time per one Bregman iterations (time/Breg (s)) in seconds. Ran on Intel Xeon CPU E5-4650 v2 @ 2.40GHz.}
\label{tab:performance}
\end{table}

\subsubsection{Qualitative tests}
\begin{figure}
  \begin{subfigure}{0.32\textwidth}
  \includegraphics[width=\textwidth]{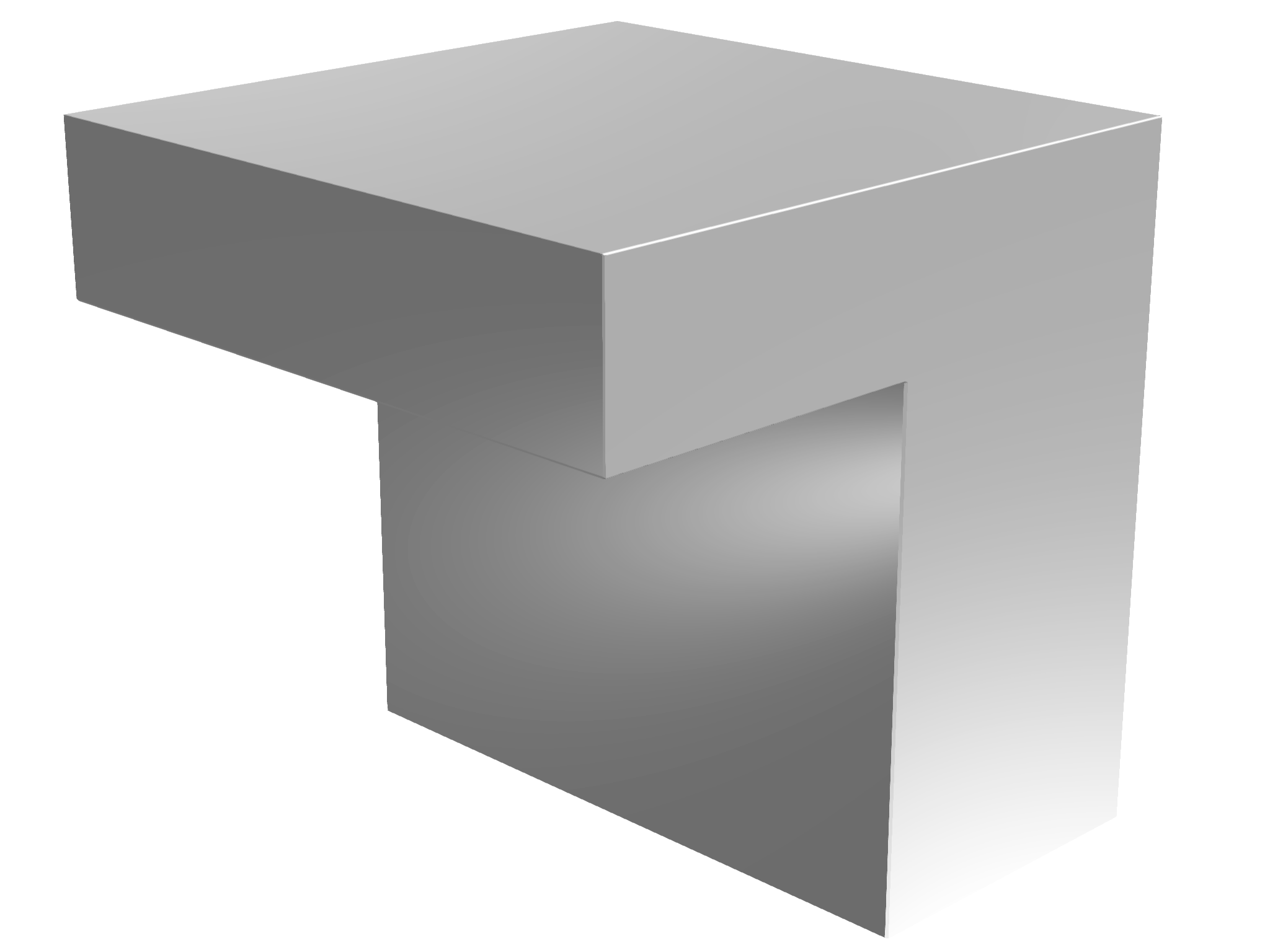}
  \end{subfigure}
  \begin{subfigure}{0.32\textwidth}
  \includegraphics[width=\textwidth]{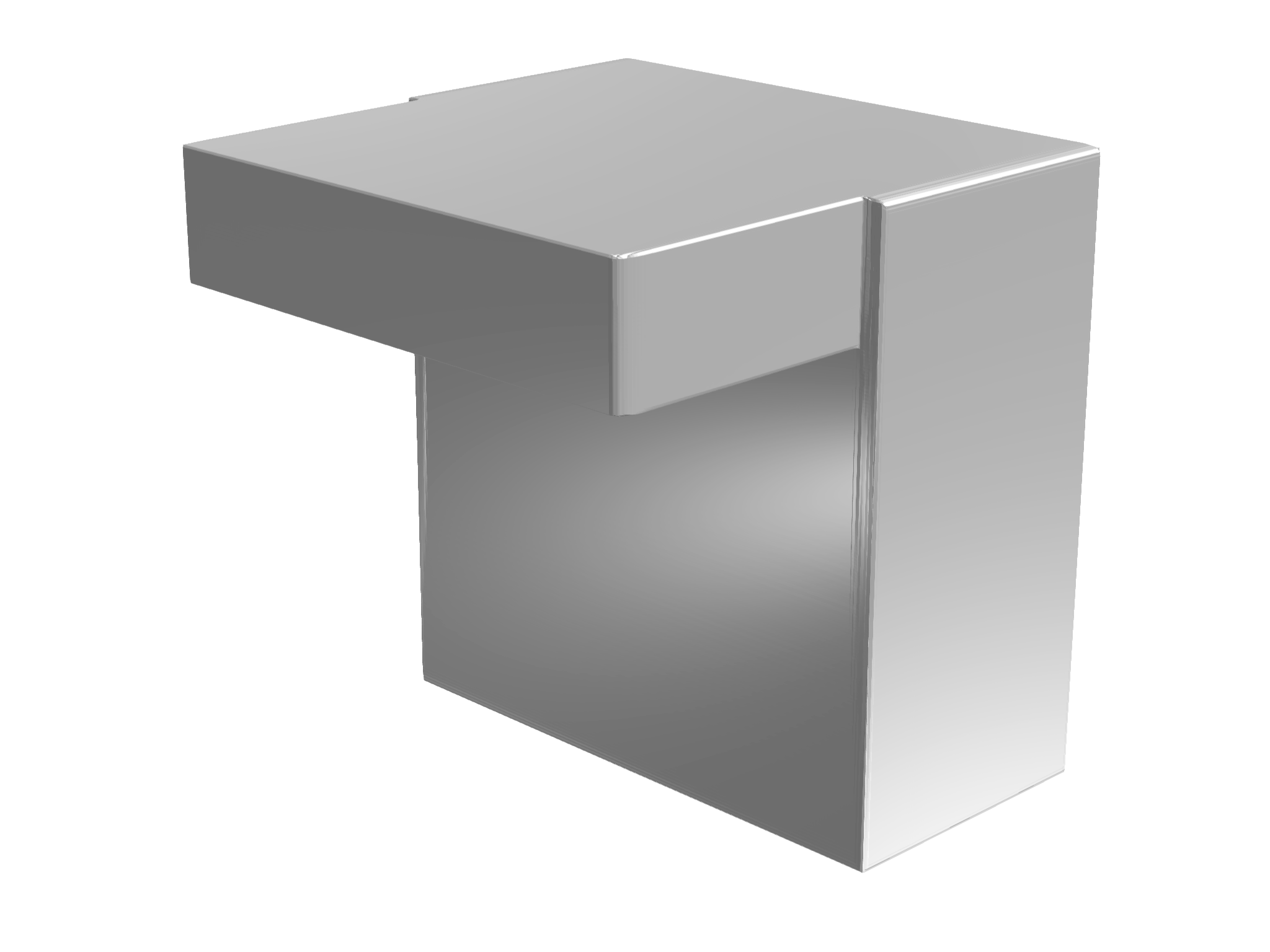}
  \end{subfigure}
  \begin{subfigure}{0.32\textwidth}
  \includegraphics[width=\textwidth]{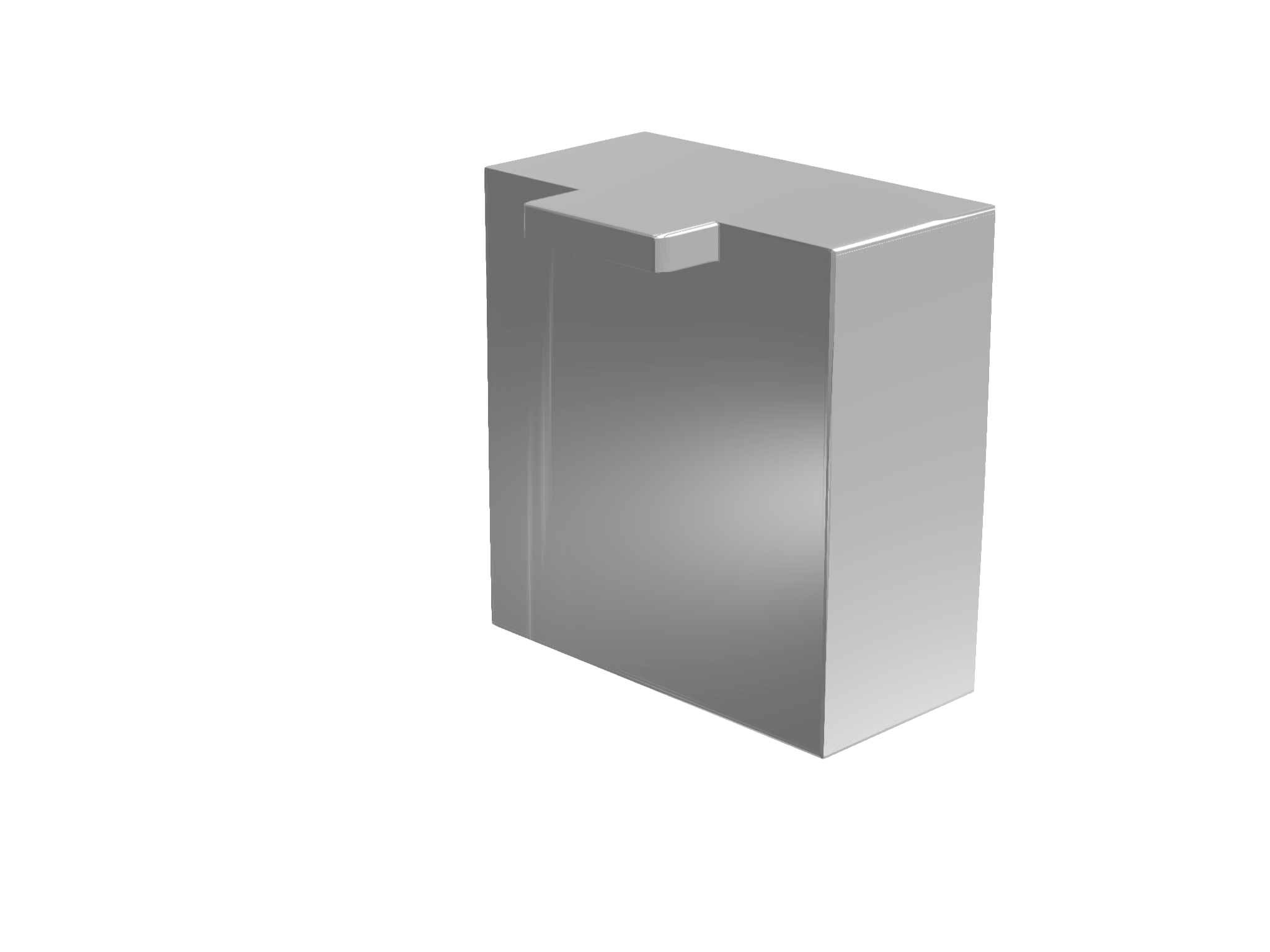}
  \end{subfigure}
  \caption{Facet breaking example in 3D \cite{BNP99} at three selected times---L-shaped facets break into
    rectangular facets.}
  \label{fig:breaking}
\end{figure}

\begin{figure}
  \begin{subfigure}{0.32\textwidth}
  \includegraphics[width=\textwidth]{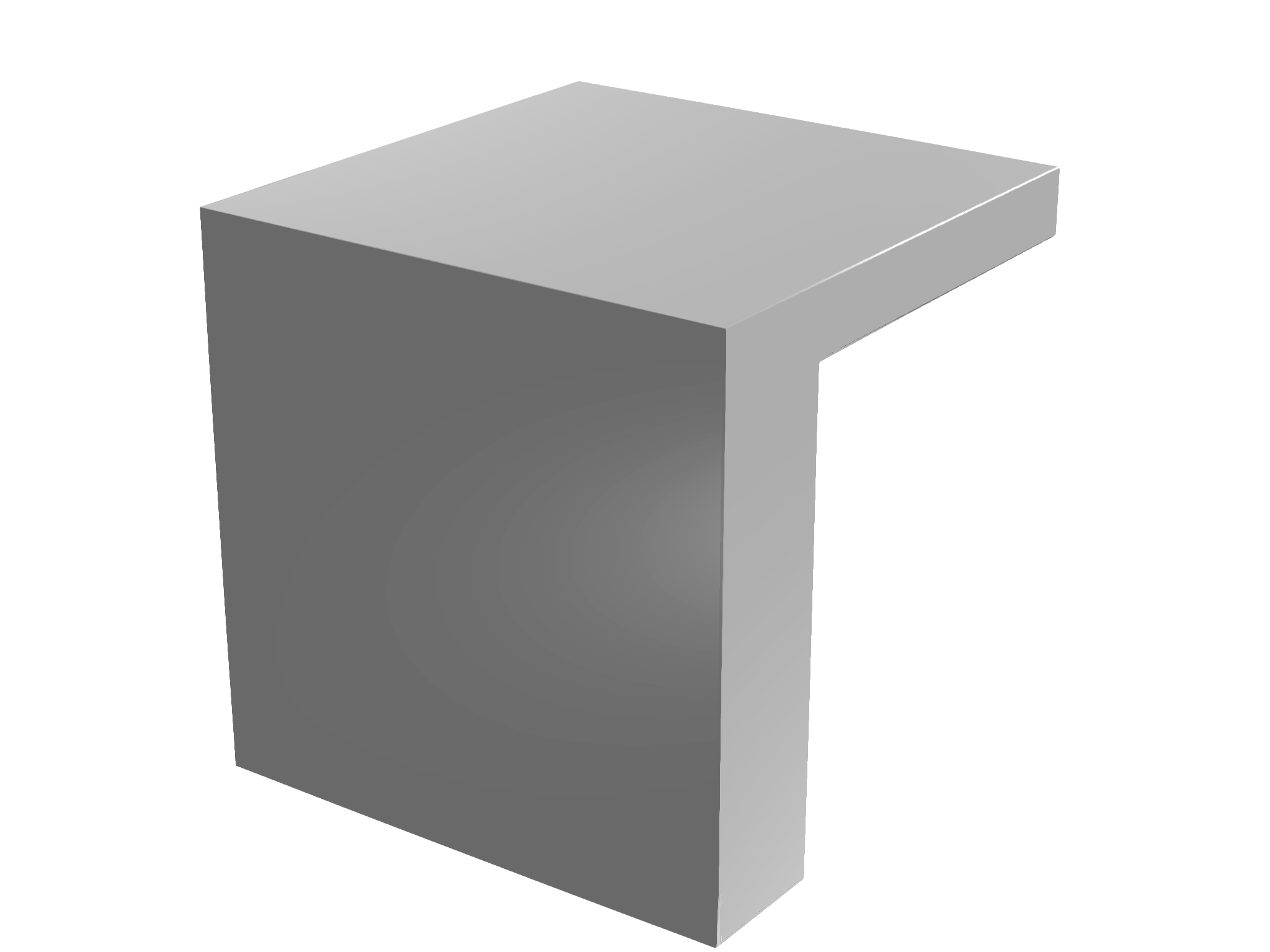}
  \end{subfigure}
  \begin{subfigure}{0.32\textwidth}
  \includegraphics[width=\textwidth]{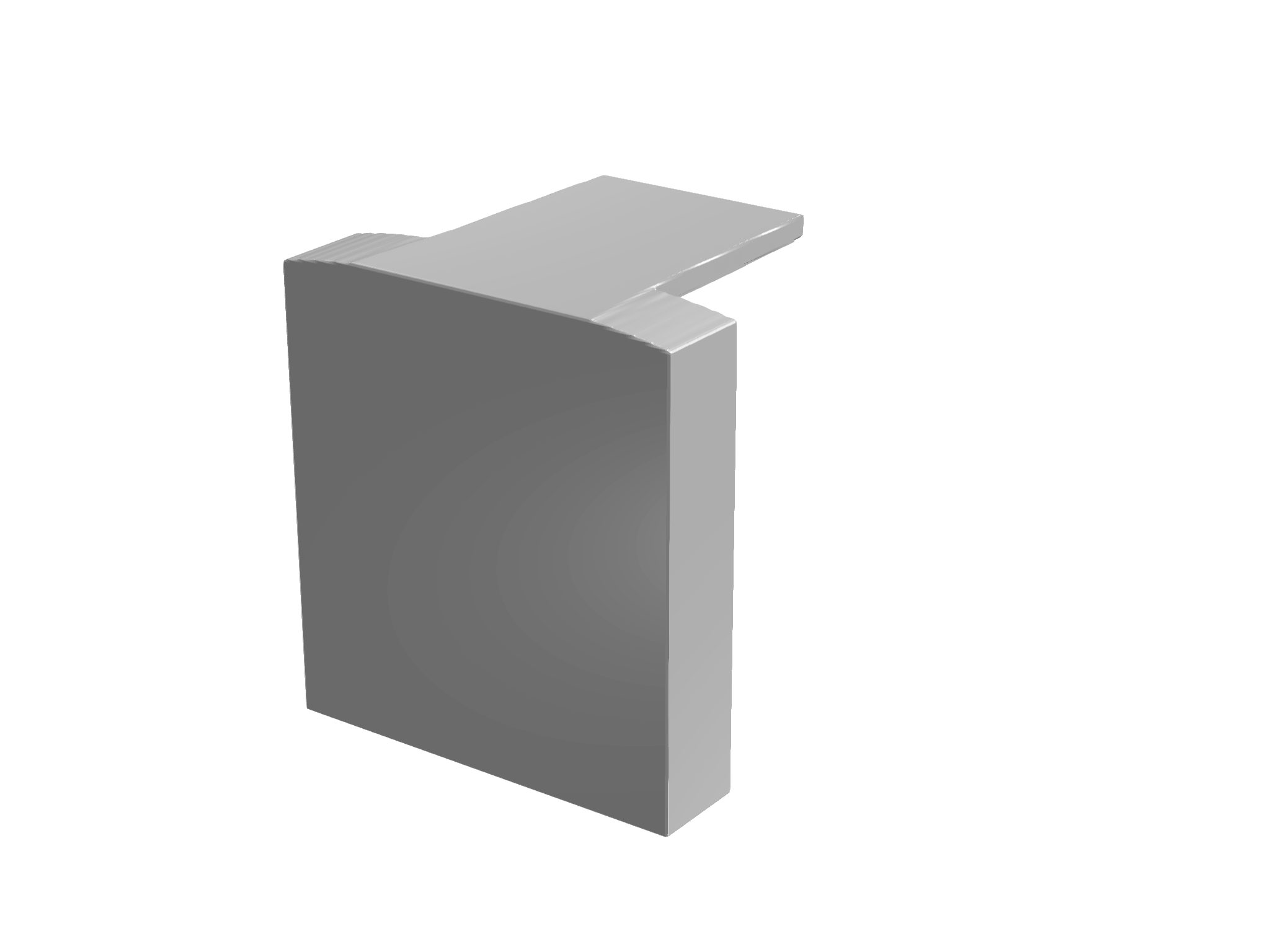}
  \end{subfigure}
  \begin{subfigure}{0.32\textwidth}
  \includegraphics[width=\textwidth]{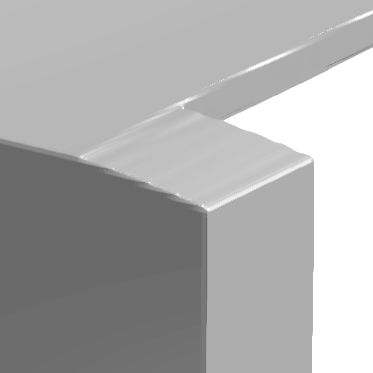}
  \end{subfigure}
  \caption{Facet bending (the top facet) in an axes-aligned initial shape for the cubic ($\ell^1$) anisotropy \cite{Lasica-comm}, with a detail on the right.}
  \label{fig:bending-L}
\end{figure}

\begin{figure}
  \begin{subfigure}{0.32\textwidth}
    \includegraphics[width=\textwidth]{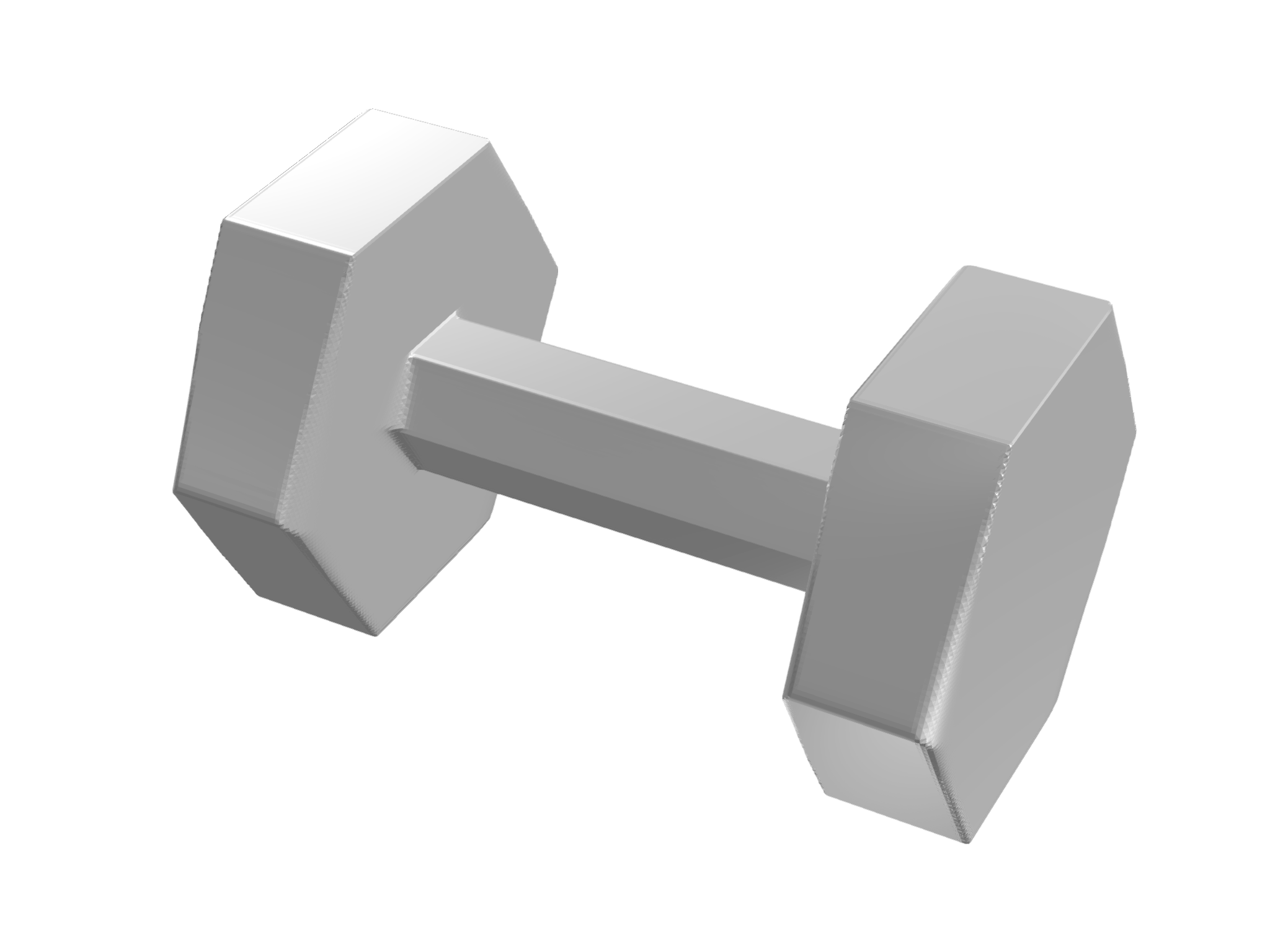}
  \end{subfigure}
  \begin{subfigure}{0.32\textwidth}
    \includegraphics[width=\textwidth]{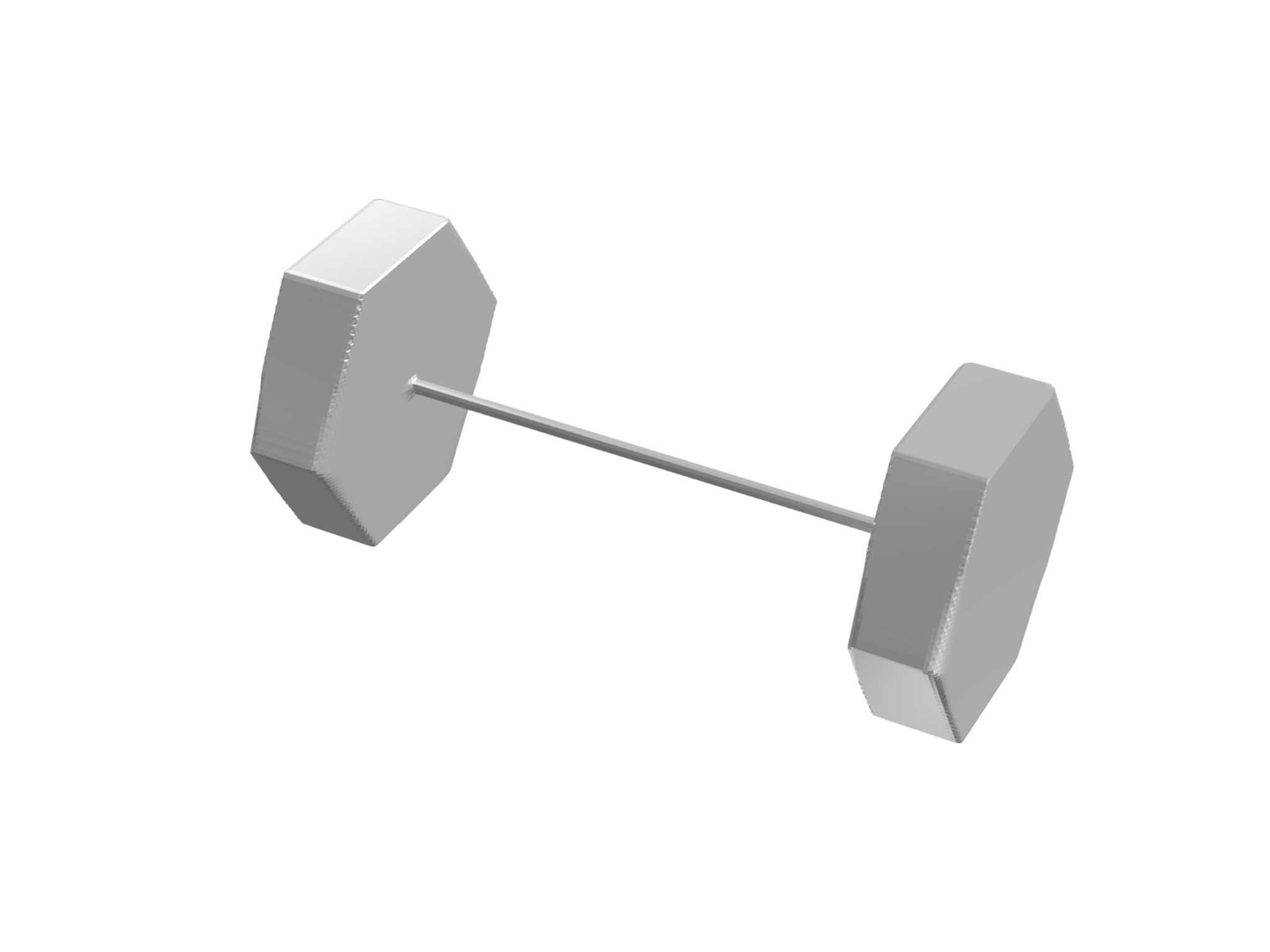}
  \end{subfigure}
  \begin{subfigure}{0.32\textwidth}
    \includegraphics[width=\textwidth]{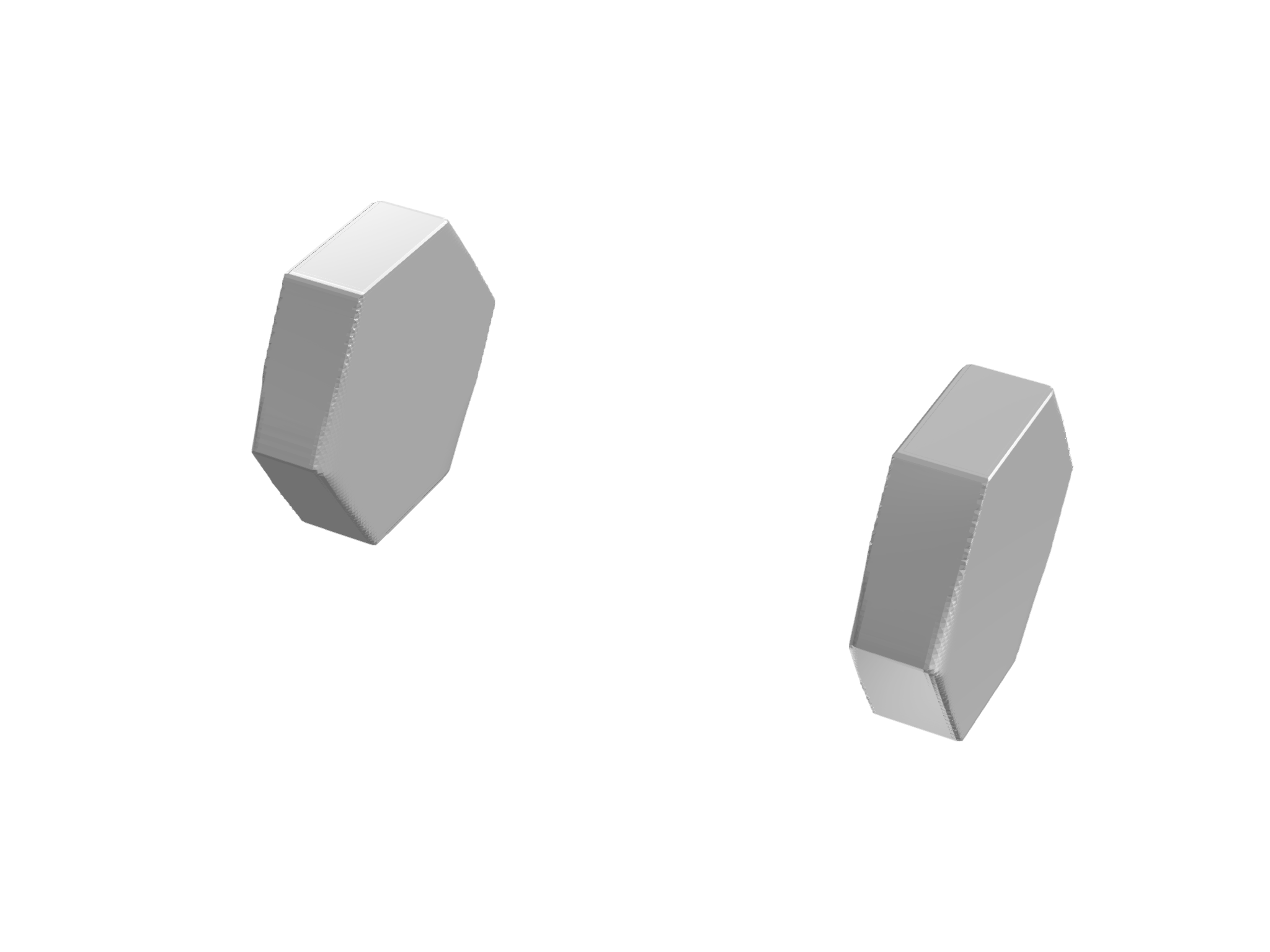}
  \end{subfigure}
  \caption{Topological change in 3D with a hexagonal anisotropy---a pinch-off of an initially
    connected dumbbell shape. FDM discretization; some rounding of edges is apparent.}
  \label{fig:dumbbell}
\end{figure}

\begin{figure}
  \includegraphics[trim={1in 1in 1in 1in},clip,width=\textwidth]{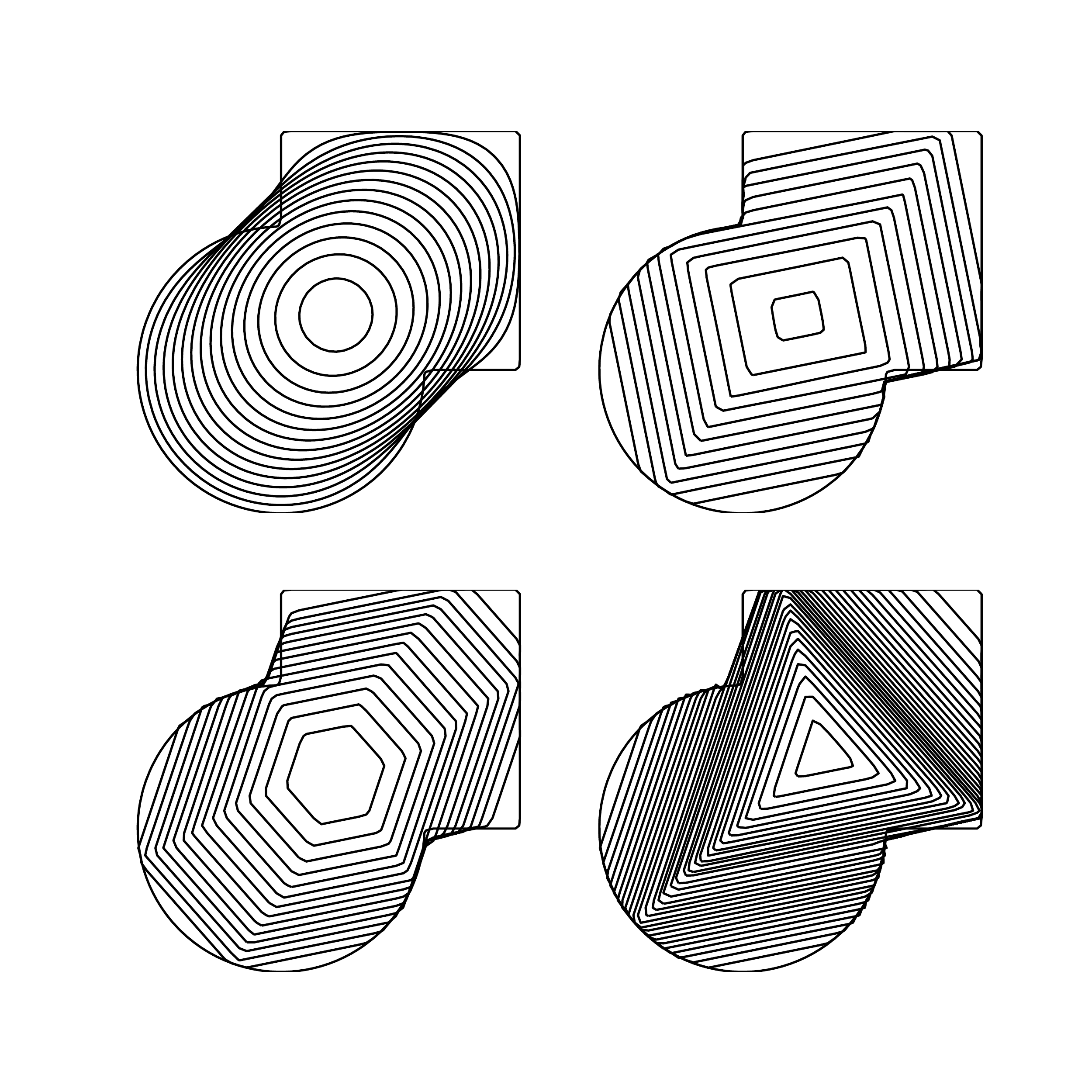}
  \caption{Flow $V = \kappa_\sigma$ with various anisotropies: $h = 10^{-4}$, $M = 64$, FEM discretization and plot
    step $0.005 = 50h$. The bounding box of the initial curve is the square $[-0.4, 0.4]^2$ inside the computational domain $(-\frac12, \frac12)^2$. Even at this relatively low resolution, artifacts caused by the redistance are not apparent.}
  \label{fig:many_squares}
\end{figure}

\begin{figure}
  \begin{subfigure}{0.32\textwidth}
    \includegraphics[width=\textwidth]{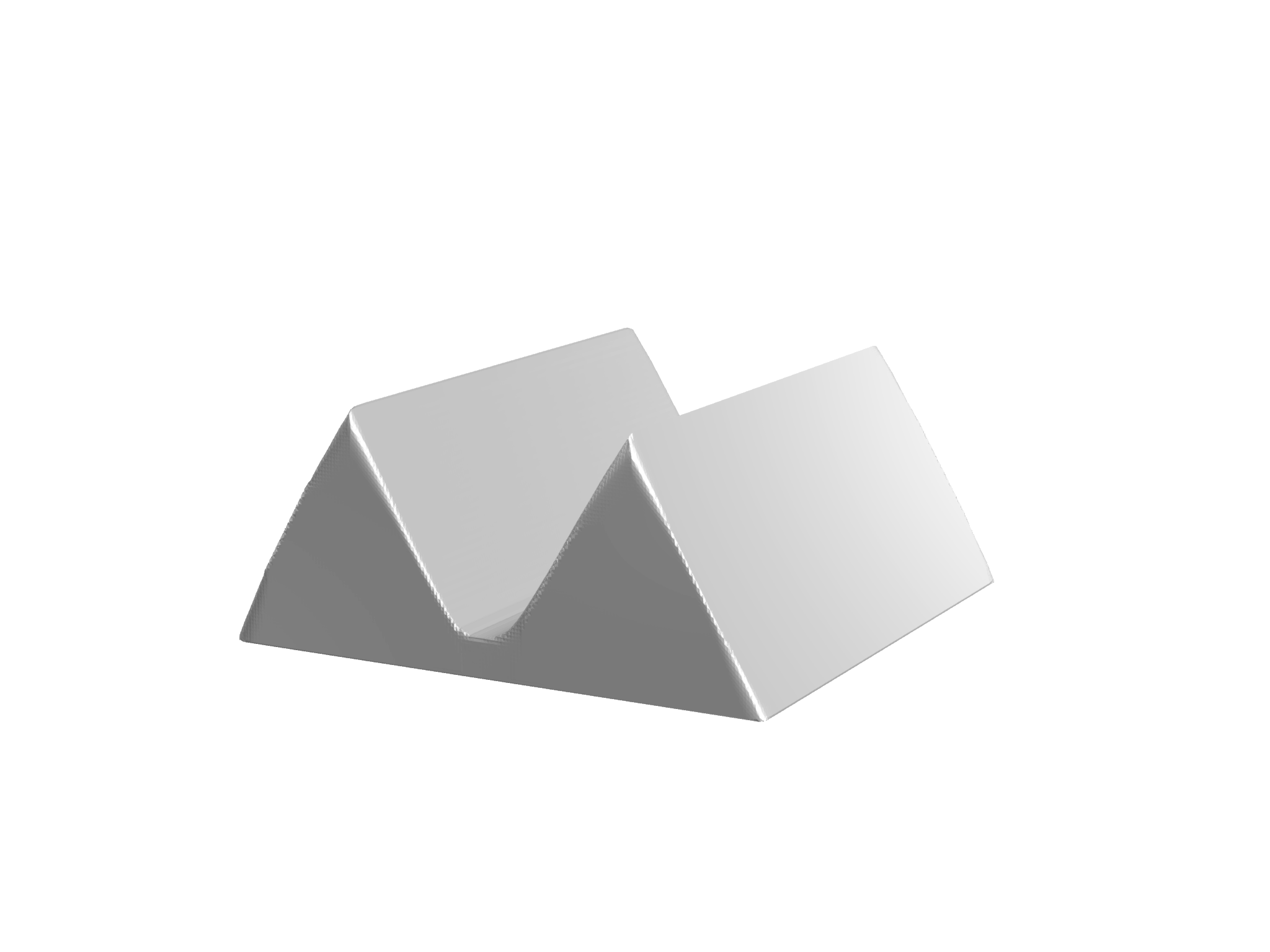}
  \end{subfigure}
  \begin{subfigure}{0.32\textwidth}
    \includegraphics[width=\textwidth]{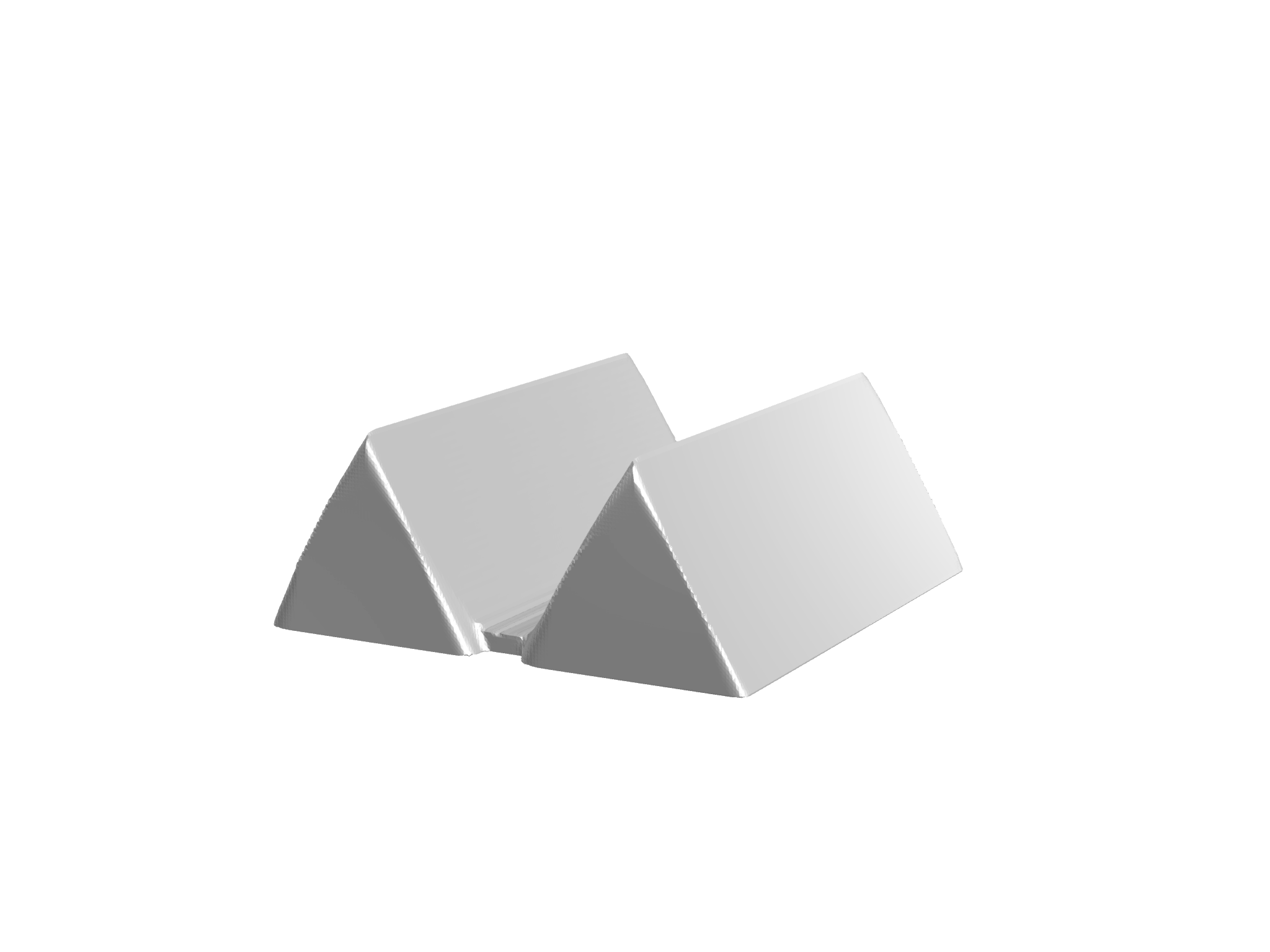}
  \end{subfigure}
  \begin{subfigure}{0.32\textwidth}
    \includegraphics[width=\textwidth]{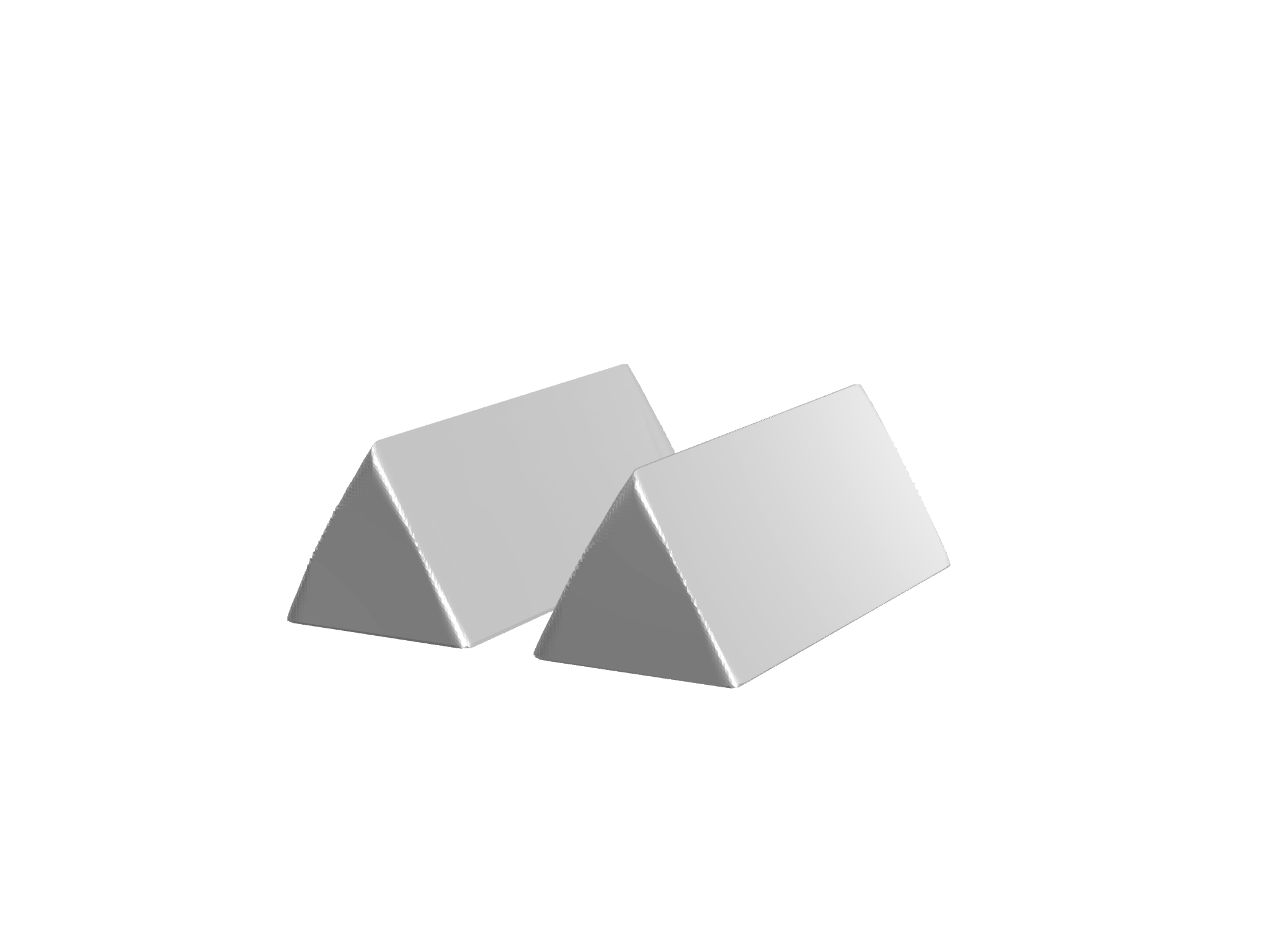}
  \end{subfigure}
  \caption{Topological change in 3D with a triangular anisotropy---a pinch-off due to a collision of
a facet and an edge. FEM discretization.}
  \label{fig:two_tris}
\end{figure}

Samples of evolutions in two dimensions for various anisotropies are shown in Figure~\ref{fig:many_squares}. In this case, the stopping condition is taken as $\norm{v_{k+1} - v_k}_{\ell^2} < 10^{-4}$ to illustrate the effect of choosing too large of a stopping tolerance, and therefore the very long edges in the bottom right figure are not completely straight. The stationary parts of the boundary are very well preserved and redistancing artifacts are not noticeable. In the rest of the figures, we show qualitative tests of facet breaking, bending and topological changes in $n = 3$, see Figure~\ref{fig:breaking}, Figure~\ref{fig:bending-L}, Figure~\ref{fig:dumbbell} and Figure~\ref{fig:two_tris}.

\subsubsection{Discussion}

The implementation of the numerical method presented in the paper is able to reproduce all of the tested features of the highly singular crystalline mean curvature flow without the need to regularize the crystalline curvature: flat faces, sharp edges and vertices, correct facet breaking and bending, and finally topological changes like neck pinching and facet-edge or edge-vertex collisions. Moreover, the method appears to be first order accurate both in space and time. We also presented a FEM discretization of the total variation energy that seems to perform visually somewhat better than the FDM discretization for non-cubic anisotropies for the cost of increased computational time. Our scheme for the reinitialization of the distance function seems to perform well without introducing redistancing artifacts.

Various optimizations, such as performing the computation only in a small neighborhood of the level
set to significantly reduce the computational complexity, as well as the coupling of the curvature flow
with the heat equation via the Gibbs-Thomson relation are under investigation \cite{PozarRIMS}.

\subsection*{Acknowledgments}
Parts of this paper are based on an extended abstract for a talk that the author gave at the 2017 Spring Meeting of the Mathematical Society of Japan.

\clearpage

\end{document}